\RequirePackage{fix-cm}
\documentclass{svjour3}                  % onecolumn (standard format)
\smartqed  % flush right qed marks, e.g. at end of proof

\usepackage{amscd}
\usepackage{amsmath}
\usepackage{amssymb}
\usepackage{graphicx}%
\usepackage{siunitx}

\usepackage{newtxtext}

\usepackage[colorlinks=true, pdfstartview=FitV, linkcolor=blue,
            citecolor=blue, urlcolor=blue]{hyperref} 
\usepackage{subcaption}
\usepackage{sidecap}
\usepackage{grffile}
\usepackage{xcolor}

\usepackage{dcolumn}% Align table columns on decimal point
\usepackage{bm}% bold math
\usepackage[mathlines]{lineno}% Enable numbering of text and display math
\usepackage{empheq}

\newcommand\at[2]{\left.#1\right|_{#2}}

\newcommand\bv[1]{\boldsymbol{#1}}

\renewcommand{\Re}{\operatorname{Re}}
\renewcommand{\Im}{\operatorname{Im}}

\newcommand{\psd}{S_{\mathrm{ESST}}}
\newcommand{\acf}{C_{\mathrm{ESST}}}

\definecolor{rred}{rgb}{0.7,0,0.1}

\definecolor{greenrb}{rgb}{0.2,0.6,0.2}

\def\d{\, \mathrm{d}}

\newtheorem{rem}{Remark}[section]
\def\br{\begin{rem}}
\def\er{\end{rem}}

\def\bi{\begin{itemize}}
\def\ei{\end{itemize}}

\numberwithin{equation}{section} 
\begin{document}

\title{Ruelle-Pollicott Resonances of Stochastic Systems in Reduced State Space. Part III: Application to the Cane-Zebiak model of the El Ni\~no-Southern Oscillation}
%\thanks{Footnote to title of article.}

% \author{Alexis Tantet}
% \affiliation{Laboratoire de M\'et\'eorologie Dynamique, Palaiseau, France}
% \email{alexis.tantet@lmd.polytechnique.fr}

% \author{Micka\"el D. Chekroun}
% \affiliation{Department of Atmospheric and Oceanic Sciences and Institute of Geophysics and Planetary Physics, University of California, Los Angeles, USA.}

% %Meteorologisches Institut, University of Hamburg, Grindelberg 5, Hamburg, Germany.}}
% %\homepage{http://www.Second.institution.edu/~Charlie.Author.}
% %\extraaffil{Institute for Marine and Atmospheric Research, Department of Physics, Utrecht University, Utrecht, The Netherlands.}
% %\address[2]{Universit\"at Hamburg, Center for Earth System Research and Sustainability, Meteorologisches Institut, Hamburg, Germany}
% \author{J. David Neelin}
% \affiliation{Department of Atmospheric and Oceanic Sciences and Institute of Geophysics and Planetary Physics, University of California, Los Angeles, USA.}
% \author{Henk A. Dijkstra}
% \affiliation{Institute for Marine and Atmospheric Research, Department of Physics, Utrecht University, Utrecht, The Netherlands.}

\author{Alexis Tantet \and Micka\"el D. Chekroun \and J. David Neelin \and Henk A. Dijkstra}
\institute{A. Tantet \at LMD/IPSL, Ecole polytechnique, Sorbonne Universite, ENS, PSL University, CNRS, Palaiseau, France\\
  \email{alexis.tantet@lmd.polytechnique.fr}
  \and
  M.D. Chekroun \and J.D. Neelin \at
  Department of Atmospheric and Oceanic Sciences and Institute of Geophysics and Planetary Physics, University of California, Los Angeles, USA
  \and
  H.A. Dijkstra \at Institute of Marine and Atmospheric Research, Department of Physics and Astronomy,  University of Utrecht, Utrecht, The Netherlands}

\date{\today}% It is always \today, today,

\maketitle

\begin{abstract}
  The response of a low-frequency mode of climate variability, El Ni\~no-Southern Oscillation, to stochastic forcing is studied in a high-dimensional model of intermediate complexity, the fully-coupled Cane-Zebiak model~\cite{Zebiak1987a}, from the spectral analysis of Markov operators governing the decay of correlations and resonances in the power spectrum.

  Noise-induced oscillations excited before a supercritical Hopf bifurcation are examined by means of complex resonances, the reduced Ruelle-Pollicott (RP) resonances, via a numerical application of the reduction approach of the first part of this contribution~\cite{Chekroun2017a} to model simulations.

  The oscillations manifest themselves as peaks in the power spectrum which are associated with RP resonances organized along parabolas, as the bifurcation is neared.
  These resonances and the associated eigenvectors are furthermore well described by the small-noise expansion formulas obtained by~\cite{gaspard2002trace} and made explicit in the second part of this contribution~\cite{Tantet2017b}.
  Beyond the bifurcation, the spectral gap between the imaginary axis and the real part of the leading resonances quantifies the diffusion of phase of the noise-induced oscillations and can be computed from the linearization of the model and from the diffusion matrix of the noise.
  In this model, the phase diffusion coefficient thus gives a measure of the predictability of oscillatory events representing ENSO.
  ENSO events being known to be locked to the seasonal cycle, these results should be extended to the non-autonomous case.

More generally, the reduction approach theorized in~\cite{Chekroun2017a}, complemented by our understanding of the spectral properties of reference systems such as the stochastic Hopf bifurcation, provides a promising methodology for the analysis of low-frequency variability in high-dimensional stochastic systems.
\end{abstract}

% % 05.10.Gg 	Stochastic analysis methods (Fokker-Planck, Langevin, etc.)
% \pacs{05.10.Gg}% PACS, the Physics and Astronomy
%                              % Classification Scheme.


% \linenumbers{}

\section{Introduction}\label{sec:introduction}

Complex and unpredictable behavior of trajectories is observed for many physical systems.
This can be due to interactions with many degrees of freedom which can be modeled by a stochastic forcing or to nonlinear coupling resulting in chaotic trajectories.
As a result, prediction beyond a certain horizon is hopeless and one focuses instead on the statistical evolution of the system.
This loss of predictability manifests itself by the evolution of an ensemble of trajectories becoming independent on its initial condition after a given time, the mixing time.
This notion of mixing in state space is in turn closely linked to the correlation function of a pair of observables, which assigns to any positive time lag the correlation between the first observable and the lagged version of the second, and thus gives a measure of the statistical dependence of the observables as the system evolves.
The power spectrum, on the other hand, describes this evolution in frequency domain.

Some chaotic or stochastic systems mix fast, in the sense that correlation functions decay exponentially with time.
As a result the corresponding power spectra are continuous.
However, systems with easily excitable modes, such as associated with persistence~\cite{tantet_early_2015}, weakly damped instabilities, or at the approach of an attractor crisis~\cite{Vaidya2008,Mauroy2016,tantet_resonances_2018,tantet_crisis_2018}, exhibit resonant behavior.
This can be seen both from the relatively slow decay of correlations after a first regime of fast decay associated with the non-resonant modes, and from peaks in the power spectra standing against a continuous background~\cite{Chekroun2014}\footnote{
  If peaks in the power spectra are continuous, the decay of correlations may be slower at first but still be exponentially fast for infinite times (as formalized by the Paley-Wiener theorem).
  In general, however, peaks may be discontinuous and prevent the exponential decay of correlations.
  Such behavior is not visible in this study.}.
In turn, resonant phenomena may be identified from such features in the correlation function and the power spectrum and eventually help better understand the physical mechanism responsible for them.
It is essential for models of these systems to resolve these resonances with appropriate time scale and frequency since the latter may explain a large fraction of the variance of the system and result in long term predictability.

The analysis of correlation functions and power spectra from observations or simulations for the study of climate variability is a common practice.
A particular important phenomenon of interest is El Ni\~no Southern Oscillation (ENSO)~\cite{Neelin1998a}, an interannual variation in the Pacific Ocean sea surface  temperature (SST).
The oscillatory nature of ENSO derives from an instability of the coupled ocean-atmosphere system, where SST anomalies cause surface-wind anomalies that in turn cause changes in surface ocean velocities.
Through the associated heat fluxes, SST is affected; in this way a positive feedback exists, the Bjerknes feedback~\cite{Jin1996b}.
The interannual time scale of ENSO is caused by an adjustment of the ocean circulation due to equatorial-wave processes, which induce a negative feedback on the SST.
The signature of ENSO is clearly visible in the periodogram estimate of the power spectrum of the Ni\~no 3.4 index~\cite[Fig.~5]{Deser2010}.
Indeed, a broad spectral peak centered around a period of 2 to 8 years stands out of a continuous spectrum.

The Cane-Zebiak (CZ) model is a system of Partial Differential Equation (PDE) with a 1.5 shallow-water ocean model coupled to a steady state atmosphere~\cite{Zebiak1987a}.
A mathematical description and analysis of the related Jin-Neelin model~\cite{Jin_al93_part1,Jin_al93_part2,Jin_al93_part3} is found in~\cite{cao2019mathematical}. 
The strength of the ocean-atmosphere interaction is controlled by a coupling parameter.
When the coupling exceeds a critical value, SST anomalies are amplified and the background climate state of the Pacific is unstable.
When the coupling strength $\xi$ exceeds a critical value $\xi_c$, a Hopf bifurcation occurs in the CZ model where a steady state loses stability to a limit cycle with a characteristic period close to the observed period of ENSO.
As such, it has been the first model to be able to simulate realistic ENSO events.
Much of the  theory developed for ENSO is based on the analysis of simulation results of this model.
Oscillations in the deterministic version of the CZ model are perfectly periodic and thus do not explain why ENSO time series are only pseudo-periodic, with a broad peak in the power spectrum.
It has been suggested by~\cite{Roulston2000} that such behavior could be seen in these models when a stochastic forcing representing fast atmospheric processes, such as westerly wind-bursts, is added.
Then noise-induced oscillations could occur before the deterministic bifurcation.
The same phenomenon has been described earlier by~\cite{wiesenfeld_effect_1982} for the van der Pol-Duffing oscillator.

In the presence of noise, new dynamical phenomena may occur, which are not explained by deterministic bifurcation theory alone.
Much can be learned about the dynamics of stochastic systems from the power spectrum or the correlation function of some observables.
For instance,~\cite{wiesenfeld_noisy_1985} study the effect of external noise on systems displaying nonlinear instabilities of periodic orbits by associating peaks in the power spectrum to Floquet exponents.
Yet, correlation functions do not give a full description of the statistical evolution of the system.
In particular, the decay of correlations strongly depends on the choice of observables.
Instead, the statistical evolution of the system is governed by the semigroup of Markov operators, which in turn can be used to compute correlation functions between any pair of observables.
In the low-dimensional case, these operators can be approximated from many short simulations or, in some cases, from a long time series.
For high-dimensional systems, this procedure is not tractable, due to the exponential increase with the number of dimensions of the number of basis functions needed to discretize the operators.
To cope with high-dimensional stochastic problems, we have introduced in~\cite{Chekroun2017a} projections of these operators on a reduced space and showed that information on the spectrum of the full Markov operators could rigorously be obtained.

In this study, we use this reduction method to analyse the spectrum of the Markov semigroup, the \emph{Ruelle-Pollicott (RP) spectrum}, for a stochastic version of the fully-coupled CZ model.
This analysis gives a new perspective on the phenomenon of noise-induced oscillations as a potential explanation of the irregularity of ENSO events.
It is found that the structure of the RP resonances undergoes a smooth yet qualitative change as the Hopf bifurcation in the CZ model is passed.
The RP resonances and the associated eigenvectors give a description of the phenomenon of noise-induced oscillations as well as of the slowing down of correlations at the approach of the bifurcation.
As the stability of the deterministic limit cycle increases, peaks in power spectra associated with resonances sharpen.
Yet, the phenomenon of phase diffusion associated with the stochastic forcing is responsible for oscillations to be irregular, or not exactly periodic.
The small-noise expansions developed in the second part of this contribution~\cite{Tantet2017b} provide a deeper understanding of this phenomenon.

In Section~\ref{sec:ergodicTheory}, we summarize the theoretical properties of the Markov semigroup, focusing on the relationship of its spectrum with the decay of correlations and the presence of resonances in power spectra.
The reduction method used to approximate this RP spectrum is also presented.
We then present in Section~\ref{sec:resultsCZ} the result of the application of this method to the CZ model.
In Section~\ref{sec:resultsJin}, we apply small-noise expansions of the RP spectrum to interpret the results obtained for the CZ model.
The results are summarized in Section~\ref{sec:conclusion}, where we also discuss the implications of our results regarding the irregularity of ENSO events.

\section{Time variability of stochastic systems and Ruelle-Pollicott resonances}\label{sec:ergodicTheory}

The authors in~\cite{Chekroun2014} have introduced a new mathematical framework to (i) understand and diagnose --- through partial observations --- the variability of turbulent flows, and (ii) to analyze parameter sensitivity that may occur in the modeling of such observations.  The framework relied on the theory of Ruelle-Pollicott (RP) resonances introduced in the mid-80's~\cite{ruelle1986locating,pollicott1986meromorphic} and that was known only by a little group of experts at the time of the publication of~\cite{Chekroun2014} working in the field dynamical system theory and the mathematical study of scattering resonances~\cite{Zworski2017}.
The RP resonances characterize the nature of the dynamics as associated with the spectrum of the underlying Liouville operator for deterministic systems or the Fokker-Planck operator for stochastic systems~\cite{gaspard2002trace,Chekroun2017a,Tantet2017b}, but are in general difficult to estimate especially if the dimension of the state space is large. 

The work~\cite{Chekroun2014} established new bridges  between the theory of RP resonances and the theory of Markov processes, once reduced state spaces are employed.
These bridges allowed for stretching new paths towards applications, especially regarding the analysis and diagnosis of complex systems' variability. 
Given an observable $h$ of a complex system and an associated reduced state space $V$ in which the observations are collected,~\cite{Chekroun2014} have shown that the eigenvalues of the dynamics' transition matrix in  $V$, called {\it reduced RP resonances}, may in fact relate to RP resonances themselves and inform about the dynamics in the full state space, once the reduced state space has been appropriately chosen.

For low-dimensional systems, Ulam's method~\cite{ulam1964collection,Dellnitz1999} is a well-known approximation of the transition matrices on a truncated set of indicator functions, which may be used to approximate the traditional RP resonances of dynamical systems.
Ulam's method led to many interesting applications in dynamical systems theory~\cite{Dellnitz1997a,Froyland2009,Koltai2010}, stochastic modeling~\cite{Froyland2013}, physical oceanography~\cite{Froyland2007,Dellnitz2009,Sebille2012,Froyland2014b} and molecular dynamics~\cite{Deuflhard1999,Schutte1999,Bittracher2015}.
Another approach, the \emph{EDMD}~\cite{Tu2013}, relies on a linear regression between snapshots of observables, making it an Extension of the Dynamic Mode Decomposition (DMD)~\cite{rowley2009spectral,schmid2010dynamic}.
This modal decomposition, also based on estimates from multi-variate time series, has been thoroughly studied by~\cite{Williams2015} who also showed applications to low-dimensional deterministic and stochastic systems.
Moreover, the link between Ulam's method and the EDMD has been stressed by~\cite{Klus2015a}, as well as the benefits and disadvantages from both methods.
Recently,~\cite{tantet_resonances_2018} applied Ulam's method to discuss the slowing down of the decay of correlations at the approach of a global attractor crisis in the Lorenz flow in terms of stable and unstable RP resonances. Together with the deterministic version of the reduction approach presented here, this work helped~\cite{tantet_crisis_2018} to analyse critical slowing down in a global attractor crisis of high-dimensional climate model.

%do not require any dictionary learning (such as for approximating the Koopman modes \cite{li2017extended}) or a sufficiently substantial amount of spatially distributed measurements \cite{tu2014dynamic}, once e.g.~an appropriate smoothing \cite{DAH} has been applied to the data prior to the estimation of RP resonances.    

% (also called an observation space)
%
%The nature of the dynamics is then characterized by the geometric organization of these eigenvalues  in the complex plane; see \cite{Chekroun2017a,Tantet2017b}.  

When the state space is high-dimensional, the reduced RP resonances, as learned within a reduced state space of low-dimension, are typically less greedy in terms of data than direct estimates in the original state space.
However, the use of a reduced state space requires more care in terms of interpretation.
In certain cases, the geometric patterns formed by reduced RP resonances in the complex plane nevertheless characterize the reduced dynamics well.
In that respect,~\cite{Tantet2017b} have shown that when the reduced RP resonances form parabolas in the complex plane, they provide an unambiguous diagnosis of an oscillation present within a stochastic or turbulent background; see also~\cite{Bagheri2014}.
 Of course for such a diagnosis to be revealed, the signal-to-noise ratio needs to be in favor of the oscillation detection, i.e.~the oscillation should not be smeared within too much of noise or small-scale turbulent behavior.
 Section~\ref{Sec_RPs} below revisits the theory of RP resonances, and Section~\ref{Sec_SL_explained0} introduces the theory of reduced RP resonances and discuss how to estimate and use them in practice to analyze dynamics through partial observations.   

 \subsection{Ruelle-Pollicott resonances and the decomposition of correlation functions and power spectra}\label{Sec_RPs}

We communicate in this section the theoretical apparatus based on the theory of {\it Ruelle-Pollicott (RP) resonances} and their numerical estimation from time series; see~\cite{Chekroun2014} for deterministic systems and~\cite{Chekroun2017a}  
for stochastic systems. We recall below the main ingredients of this theory, in particular regarding the fundamental role 
that RP resonances play in the decomposition of power spectral density (PSD) and correlation functions.

The Cane-Zebiak (CZ) model~\cite{Zebiak1987a} subject to stochastic perturbations is considered in this study; the addition of noise appearing in the wind stress.
In the deterministic setting, it corresponds to a spectral 
discretization of a system of Partial Differential Equations (PDEs).
Assuming that time is continuous, i.e.~before applying a time discretization scheme to the model  we are thus led to Stochastic Differential Equations (SDEs) of the form:
\begin{equation}\label{Eq_SDE}
\d X= F(X) \d t + G(X) \d W_t, \qquad X \in \mathbb{R}^p.%\mbox{ with } Q(x)=G(x) G(x)^T.
\end{equation}
Here $W_t=(W_t^1,\cdots,W_t^q)$ is an $\mathbb{R}^q$-valued Wiener process ($q$ not necessarily equal to $p$) whose components are mutually independent standard Brownian motions.

In Eq.~\eqref{Eq_SDE}, the drift part is provided by a (possibly nonlinear) vector field $F$ of $\mathbb{R}^p$ which, in this study, represents the CZ model.
The stochastic diffusion in its It\^o version, given by $G(X) \d W_t $, has its $i^{th}$-component ($1\leq i\leq p$) given by

\begin{equation}
\sum_{j=1}^q G_{ij} (X) \d W_t^j, \; \; q\geq 1.
\end{equation}
It corresponds in our case to the stochastic wind forcing.
In what follows we assume that the vector field $F$ and the matrix-valued function 

\begin{equation*}
G:\mathbb{R}^p \rightarrow \mbox{Mat}_{\mathbb{R}}(p\times q), 
\end{equation*}
satisfy regularity conditions that guarantee the existence and the uniqueness of mild solutions, as well as the continuity of the trajectories; see, e.g.~\cite{cerrai2001second,flandoli2010flow} for such conditions in the case of locally Lipschitz coefficients.

It is well-known that the evolution of the probability density of the random $\mathbb{R}^p$-valued variable, $X_t$ (solving Eq.~\eqref{Eq_SDE}),  is governed by the {\it Fokker-Planck equation}

\begin{equation}\label{Eq_FKP}
\partial_t  \rho(X,t) =\mathcal{A} \rho(X,t)= -\mbox{div}(\rho(X,t)F(X)) +\frac{1}{2} \mbox{div} \nabla (\mathbf{\Sigma}(X)  \rho(X,t)), \;\; X\in \mathbb{R}^p,
\end{equation}
with $\mathbf{\Sigma}(X)=G(X) G(X)^T$.
In weather forecasting and climate modelling where the use of ensemble simulations is common, one may think of the Fokker-Planck equation to govern the evolution of ensembles represented by probability densities.

What is less-known is that the spectral properties of the 2nd-order differential operator, $\mathcal{A}$, informs about  fundamental objects such as the power spectra or correlation functions computed typically along a stochastic path of  Eq.~\eqref{Eq_SDE}.
We briefly recall the main elements hereafter and refer to~\cite{Chekroun2014,Chekroun2017a} for more details. 

First, given an observable $h:\mathbb{R}^p\rightarrow \mathbb{R}$ for the system $\eqref{Eq_SDE}$, we recall that  the {\it power spectrum} $S_h(f)$ is obtained by taking the Fourier transform of the correlation function $C_h (t)$, namely 

\begin{equation}\label{Eq_corr}
S_h(f) =\widehat{C_h (t)}, \; \mbox{ with } \; C_h (t)=\lim_{T\rightarrow \infty} \frac{1}{T} \int_0^T h(X_{s+t}) h(X_{s}) \d s - \left(\lim_{T \rightarrow \infty} \frac{1}{T} \int_0^T h(X_{s}) \d s \right)^2,  
\end{equation}
where $X_t$ is a stochastic process solving~\eqref{Eq_SDE}.
Relying on the pointwise ergodic theorem, the correlation function is defined here in terms of long time averages, which, in statistics, are estimated with sample means or sample (Pearson) correlations.
% and satisfying $X_t^x=x$ at time $t=0$. 

As shown in \cite{Chekroun2017a}, for a broad class of SDEs that possess an ergodic probability distribution $\mu$,   
the spectrum, $\sigma(\mathcal{A})$, of the Fokker-Planck operator, $\mathcal{A}$,  is typically contained in the left-half complex plane, $\{z\in \mathbb{C}\;:\: \Re(z)\leq 0\}$ and its resolvent  $R(z)=(z \mbox{Id}-\mathcal{A})^{-1}$, is a well-defined linear operator that satisfies

\begin{equation}\label{Eq_PSD_RPs}
S_h(f) =\int_{\mathbb{R}^p} h(X) \big[R(i f) h \big](X)\d \mu.
\end{equation}
In practice, the functional space in which $h$ in~\eqref{Eq_corr} or~\eqref{Eq_PSD_RPs} lies, and in  which the eigenfunctions of $ \mathcal{A}$, are considered plays a key role to rigorously define the spectrum $\sigma(\mathcal{A})$.
We refer to~\cite{Chekroun2017a} for the reader interested by these aspects.

\begin{figure}
  \centering
  \includegraphics[width=.4\textwidth]{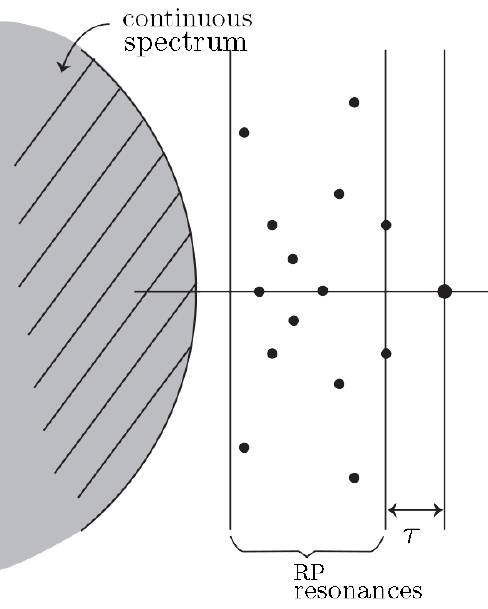}
  \caption{
    Schematic of the spectrum of $\mathcal{A}$ given in~\eqref{Eq_FKP}.
    The Ruelle-Pollicott (RP)  resonances are the isolated eigenvalues of the Fokker-Planck operator, $\mathcal{A}$, defined in~\eqref{Eq_FKP}; they are represented by red dots in Panel (b) and by black dots here.
    The rightmost vertical line represents the imaginary axis above which the power spectrum  lies; see Panel (a) for another perspective.
    The rate of decay of correlations is controlled by the spectral gap $\tau$; see~\cite{froyland1997computer,Chekroun2014}. [Courtesy of Maciej Zworski.]}\label{RP_schema}
\end{figure}

The resolvent can be understood as the Laplace transform of the time evolution generated by the Fokker-Planck operator.
Its poles thus relate to the stability of the solutions of the Fokker-Planck equation.
In~\eqref{Eq_PSD_RPs}, the frequency $f$ lies in the complex plane $\mathbb{C}$, and the poles of the resolvent  $R(i f)$ \textemdash which correspond to the RP resonances \textemdash introduce singularities into $S_h(f)$.
Once the power spectral density (PSD) is calculated, i.e.~once $|S_h(f)|$ is computed with $f$ taken to be real, these poles manifest themselves as peaks that stand out over a continuous background at the frequency $f$ if the corresponding RP resonances with imaginary part $f$ (or nearby) are close enough to the imaginary axis.
The RP resonances are the isolated eigenvalues of finite multiplicity of $\mathcal{A}$;
Although this is not the case in the applications considered here, 
the spectrum of the infinite-dimensional operator $\mathcal{A}$ may contain a continuous part lying typically in a sector $\{z\in \mathbb{C}\;:\: \Re(z)\leq -\gamma\}$ (for some $\gamma >0$); see Fig.~\ref{RP_schema}.

Denoting by $\lambda_j$s, the $N_p$ poles of the resolvent $R(z)$ of $\mathcal{A}$, i.e.~the RP resonances, and by 
$m_1,\cdots,m_N$ their corresponding orders, the correlation function, $C_h (t)$, possesses then the following expansion

\begin{equation}\label{Eq_decomp_corr1}
C_{h}(t)= \sum_{j=1}^{N_p} \Big[\sum_{k=0}^{m_j-1} \frac{\tau^k}{k!} (\mathcal{A}-\lambda_j \textrm{Id})^{k}\Big] a_j(h) e^{\lambda_j t} + \mathcal{Q}(t),
\end{equation}
where $\mathcal{Q}(t)$ exhibits typically a decay property associated with properties of the essential spectrum of $\mathcal{A}$ (see~\cite{Engel2001}, Chap.~IV.1, for a definition).
Even if the $\lambda_j$'s do not depend on the observable $h$, the 
coefficients, $a_j(h)$'s, in the expansion~\eqref{Eq_decomp_corr1}, do; see~\cite[Corollary 2.1]{Chekroun2017a}. More precisely, denoting by $\psi_j$ the eigenfunctions associated with the  $\lambda_j$'s (and by $\psi_j^\ast$ those associated with the adjoint of $\mathcal{A}$), we have

\begin{align}
	a_j(h) = \langle \psi_j^*, h \rangle_\mu \langle h, \psi_j\rangle_\mu,
	\label{eq:spectralWeights}
\end{align}
where  $\langle f, g\rangle_\mu=\int_\mathcal{H}f(x)  g(x) \mu(\d x)$, in which $\mu$ is the ergodic statistical equilibrium associated with Eq.~\eqref{Eq_SDE}.
In other words, $\mu$ is the measure giving the long-term averages of observables.
For conditions ensuring the existence of such a statistical equilibrium  we refer to~\cite[Sec.~2]{Chekroun2017a}. 

If we assume that $\Re(\lambda_j) < 0$ for $j > 0$, each $\lambda_j$ is simple and the absence of an essential spectrum for $\mathcal{A}$, then 
the correlation $C_{h}(\tau)$ in~\eqref{Eq_decomp_corr1} takes the simpler form of a weighted sum of complex exponentials, and the corresponding power spectrum $S_{h}(f)$ possesses itself a similar decomposition in terms of Lorentzian functions, namely:

\begin{align}
	S_{h}(f) = - \frac{1}{\pi} \sum_{j = 1}^\infty a_j(h) \frac{\Re(\lambda_j)}{(f - \Im(\lambda_j))^2 + \Re(\lambda_j)^2}, \qquad f\in \mathbb{C}. 
	\label{eq:spectralPower}
\end{align}

As mentioned earlier, the RP resonances have been introduced by Ruelle~\cite{ruelle1986locating} and Pollicott~\cite{pollicott1986meromorphic} for discrete and continuous chaotic deterministic systems; see also~\cite{Butterley_Liverani} for the case of Anosov flows.
\cite{Chekroun2017a} have recently shown that the theory of RP resonances extends to a stochastic framework giving access to stochastic-analysis tools to justify decomposition formulas such as~\eqref{Eq_decomp_corr1} or~\eqref{eq:spectralPower} for a broad class of SDEs.
These results are based on the theory of Markov semigroups on one hand, and the spectral theory of semigroups, on the other; see also~\cite{gaspard2002trace,dyatlov2015stochastic} for complementary approaches.

All eigen-spectra considered here and in part II of this contribution are discrete and of finite multiplicity, i.e.~composed of RP resonances only.
However, not all SDEs have a resonances (see e.g.~\cite{avram_spectral_2013}) and not all SDEs have an essential spectrum (e.g.~Ornstein-Uhlenbeck processes,~\cite{Pavliotis2014}, Chap. 4).
The resonances of the spectrum is associated with a (fast) exponential decay of correlations.
The essential part, associated with the term $\mathcal{Q}(t)$ in~\eqref{Eq_decomp_corr1}, may instead be associated with an algebraic decay of correlations.
\begin{figure}
	\centering
	\includegraphics[width=.6\textwidth]{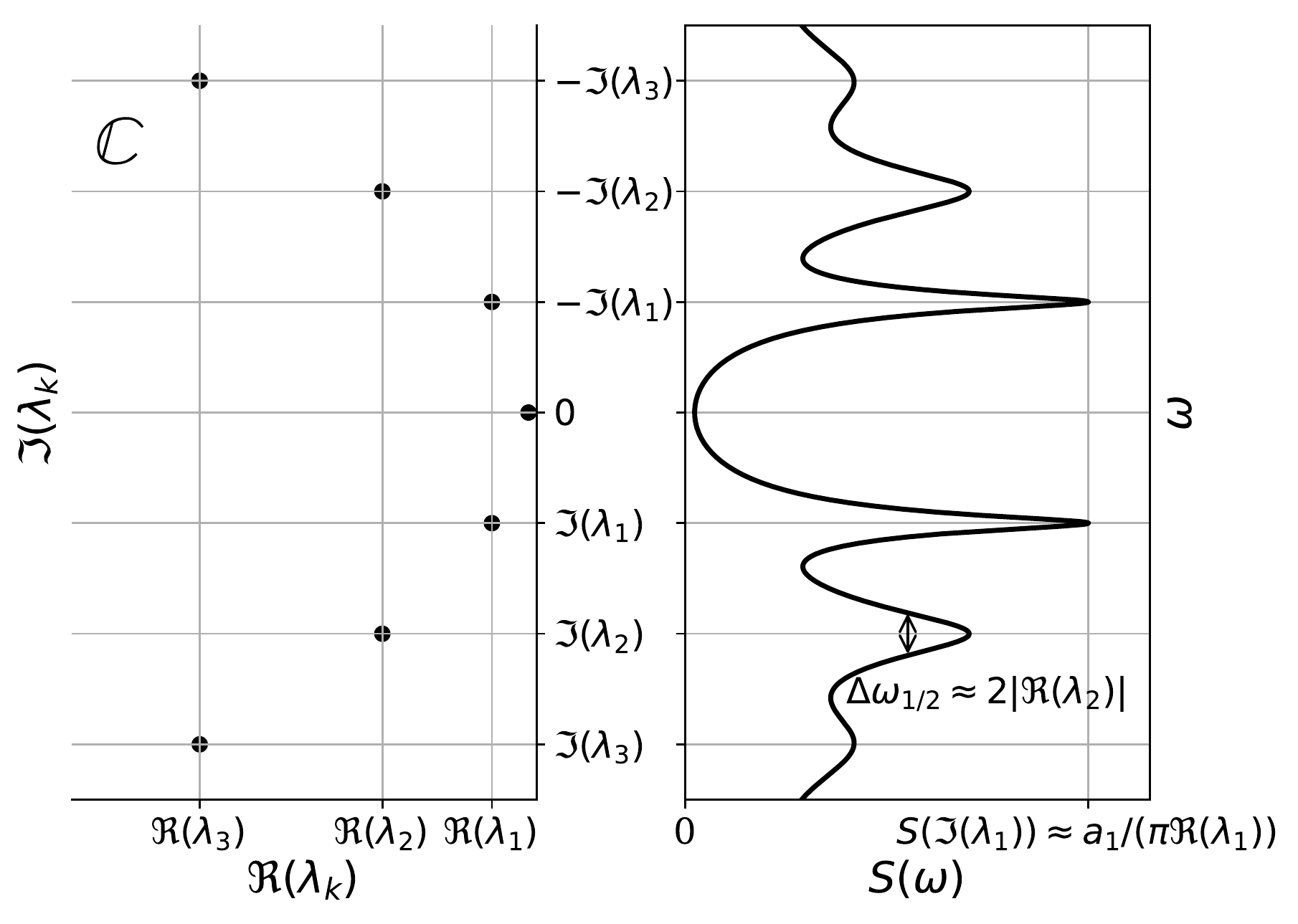}
	\caption{Schematic of the PSD as a combination of Lorentzian profiles associated with each RP resonance (pole) in~\eqref{eq:spectralPower}.}\label{fig:schemPower}
\end{figure}

In summary, the decompositions~\eqref{Eq_decomp_corr1} and~\eqref{eq:spectralPower}
combined with Fig.~\ref{fig:schemPower} inform us about the following features:
%%%%%%%
\begin{itemize}
	\item[(i)] Each RP resonance is associated with an exponential contribution to the decay of correlation.
	\item[(ii)] The closer an eigenvalue to imaginary axis, the slower the decay.
	\item[(iii)] In the limit of purely imaginary eigenvalues, the associated contributions
	to the correlation functions are purely oscillatory and prevent the decay of correlations.
	\item[(iv)] The angular frequency at which each contribution oscillates
	is given by the imaginary part of the associated eigenvalue.
	\item[(v)] Eigenvalues close to the imaginary axis are associated with resonances (i.e.~peaks) in the
	power spectrum.
	The spectral peak is located at the frequency given by the imaginary part of the eigenvalue
	and its width is proportional to the absolute value of the real part.
	\item[(vi)] The contribution of each eigenvalue to a correlation function or a power spectrum
	is weighted as in~\eqref{eq:spectralWeights}, corresponding to the projection
	of the observables $h$ onto the eigenfunctions of $\mathcal{A}$ and its adjoint.
\end{itemize}

\subsection{Estimating resonances and diagnosing the reduced dynamics from time series: The reduced RP resonances}\label{Sec_SL_explained0}

Because they are related to the Fokker-Planck operator, $\mathcal{A}$, the RP resonances inform on the structure of the underlying SDE. A natural question then arises: is it possible to infer the ``shape'' of an SDE (i.e.~$F$ and $G$ in   \eqref{Eq_SDE}) from the observation of its solutions? The short answer to this question is ``No" in general, except in certain cases (see e.g.~\cite{crommelin2006reconstruction}), but as we will see, a reduced spectrum associated with the observation space and related to RP resonances, can be estimated from time series.  Depending on the geometric patterns that it forms in the complex plane, this spectrum may inform us in turn about the sought structural ingredients, in terms of reduced dynamics.  In that respect, the approach of \cite{Chekroun2014} paved a roadmap for addressing this question from a practical viewpoint, while \cite{Chekroun2017a} analyzed in greater details its theoretical foundations, and \cite{Tantet2017b} provided general understanding regarding the type of geometric patterns that RP resonances may exhibit in the case of the stochastic Hopf bifurcation.

In practice, only partial observations of the solutions to~\eqref{Eq_SDE} are available, e.g.~few solution's components. Speaking roughly,~\cite[Theorem A]{Chekroun2014} shows that from partial observations of a complex system  that lie in a reduced state space $V$ and are taken at a sampling rate $\tau$, a (reduced) Markov operator $\mathfrak{T}_\tau$ with state space $V$ can be inferred from these observations such that (i) $\mathfrak{T}_\tau$ characterizes the coarse-graining in $V$ of the 
transition probabilities in the full state space, and (ii) the spectrum of $\mathfrak{T}_\tau$ relates to the RP resonances, at best in an averaged sense; see also~\cite[Theorems 3.1 and 3.2]{Chekroun2017a} and~\cite{schutte1999direct,froyland2014computational}.
The operator $\mathfrak{T}_\tau$ can be thought as the operator that maximizes the likelihood of the evolution during a time $\tau$ of densities in the reduced state space conditionned on information in the reduced phase space only.
In practice the dimension of $V$ is kept low so that $\mathfrak{T}_\tau$ can be efficiently estimated from time series via a maximum likelihood estimator (MLE). 
The reduced state space $V$ should be also chosen such that the observed dynamics in $V$ carry relevant information on the low-frequency variability of the system.
For instance, in the case of the CZ model analyzed in Sec.~\ref{sec:resultsCZ}, we define a two-dimensional observable, $h$,  from components given by the eastern sea surface temperature, on one hand, and the western thermocline depth, on the other.

We detail below our estimation procedure and what (ii) means. 
First a bounded domain $\mathcal{D}$ of $V$ should be chosen large enough so that ``most realizations'' of the stochastic process $X_t$ solving Eq.~\eqref{Eq_SDE} fall inside $\mathcal{D}$ after application of the observable $h: \mathbb{R}^p \rightarrow V$,  i.e. $\mathcal{D}$ must be chosen so that $h(X_t)$ belongs to $\mathcal{D}$ for many realizations of the noise in Eq.~\eqref{Eq_SDE}.  This domain is then discretized as the union of $M$ disjoint boxes $B_j$, forming thus a partition. 

%The resulting grid should be fine enough to resolve the
%leading eigenvectors. 

In practice our observations are made out of the process $X_t$ at discrete time instants
$t=t_n$,  given as multiple of a sampling time $\tau$, i.e.~$t_n=n\tau$ with $1\leq n\leq N$, with $N$ assumed to be large. 
These observations made in the observation space $V$ are denoted by $Y_n=h(X_{t_n})$.
As explained in the Supporting Information of~\cite{Chekroun2014}, the Markov operator $\mathfrak{T}_\tau$ is approximated by the $M\times M$ transition matrix $T_\tau$ whose entries are given by the transition-probability estimates

\begin{equation}\label{P_estimator}
(T_{\tau})_{ij}=\frac{\#\bigg\{\Big( Y_{n}\in  B_j \Big)\wedge \Big(Y_{n+1} \in  B_i \Big)\bigg\}}{\#\Big\{Y_{n} \in  B_j \Big\}},
\end{equation}
%%%
where the $B_j$'s form a partition (composed of $M$ disjoint boxes) of the aforementioned domain $\mathcal{D}$ in $V$; see e.g.~\cite{schutte1999direct,crommelin2006fitting,Chekroun2014}. In~\eqref{P_estimator}, the notation $\#\{(Y_{n}\in  B_k)\}$ gives the number of observations $Y_{n}$ falling in the box $B_k$, and the logical symbol ``$\wedge$'' means ``and." The leading eigenvalues of the transition matrix $T_\tau$ can then be computed
with an iterative algorithm such as ARPACK~\cite{Lehoucq1997}.

In practice, we are interested with the eigenvalues $\lambda_k(\tau)$ obtained 
from the eigenvalues $\zeta_k(\tau)$ of the Markov matrix $T_\tau$, according to

\begin{equation}\label{Eq_lambda}
  \lambda_k(\tau)
  = \frac{1}{\tau} \log \zeta_k(\tau)
  = \frac{1}{\tau} \left(\log \left| \zeta_k(\tau)\right| 
  +  i\arg \zeta_k(\tau) \right), \qquad 1\leq k\leq M,
\end{equation}
where $\arg(z)$ (resp.~$\log(z)$) is the principal part of the argument (that we adopt to lie in $[-\pi, \pi)$ in this article) (resp.~logarithm) of the complex number $z$.
At a basic level, the motivation behind~\eqref{Eq_lambda} is that the eigenvalues of $T_\tau$ as the eigenvalues of a Markov matrix, lie within the unit circle (representation that was adopted in~\cite{Chekroun2014}) whereas we want here to relate these eigenvalues with the RP resonances associated with the original Fokker-Planck operator $\mathcal{A}$.  Thus, the $\lambda_k(\tau)$'s given by~\eqref{Eq_lambda} lying naturally in the left half of the complex plane, offer a direct comparison with the RP resonances. 

The eigenvalues $\lambda_k(\tau)$ are called hereafter the {\it reduced RP resonances}. The precise relationships between the $\lambda_k(\tau)$'s and the actual RP resonances are non-trivial to characterize in general. Nevertheless, in certain cases, the reduced RP resonances are very useful to diagnose and characterize important dynamical features such as nonlinear oscillations embedded within a stochastic background.   

Indeed the rigorous result~\cite[Theorems 3.1]{Chekroun2017a} ensures that $T_\tau$  characterizes entirely the coarse-grained probability transitions --- imposed by the observation space $V$ --- of the actual dynamics, and thus at an intuitive level if a dominant recurrent behavior occurs within an irregular background, then $T_\tau$ must still ``feel'' this recurrent behavior within $V$, in case this dominant behavior is reflected in $V$.  As pointed out already in~\cite{Chekroun2014} such a recurrent behavior is manifested by eigenvalues of $T_\tau$ distributed evenly along an inner circle typically close to the unit circle, or forming a parabola in the complex plane if one look at the $\lambda_k(\tau)'s$ as explained in details in~\cite{Tantet2017b}.

Now the fact that $T_\tau$  characterizes the coarse-grained  probability transitions and not the actual probability transitions must be 
an important feature to keep in mind. For instance if important dynamical features occur in part of the state space not included in $V$, then the reduced RP resonances will not accurately reflect the actual RP resonances.
%; see \cite[Theorem 3.2]{Chekroun2017a}

In fact, and as pointed out in~\cite[Sec.~3]{Chekroun2017a}, even if the family of Markov matrices  $(T_\tau)_{\tau \geq 0}$ satisfies the semigroup property (i.e.~$T_{\tau+\tau'}=T_\tau T_{\tau'}$), the reduced RP resonances do not match necessarily with the RP resonances.  
Nevertheless, in this case, the reduced RP resonances relate to an averaged Fokker-Planck operator.
Indeed, in the limit $N, M\rightarrow \infty$, the family of Markov operator $(\mathfrak{T}_\tau)_{\tau \geq 0}$ satisfies $\mathfrak{T}_{\tau+\tau'}=\mathfrak{T}_\tau \mathfrak{T}_{\tau'}$.
Thanks to~\cite[Theorem 3.2]{Chekroun2017a}, this property can be used to infer that the reduced RP resonances (independent of $\tau$ in this case) approximate  the point spectrum $\sigma_p(\overline{\mathcal{A}})$ of the following averaged Fokker-Planck operator, given formally for all $\Psi$ sufficiently smooth by 

\begin{equation}\label{Eq_average_FKP_ope}
\overline{\mathcal{A}} \Psi (y)=\int_{x\in h^{-1}(\{y\})} \bigg(-\mbox{div}(\Psi(y) F(x)) +\frac{1}{2} \mbox{div} \nabla \big(\mathbf{\Sigma}(x) \Psi (y)\big)  \bigg) \d \mu_y(x),
\end{equation} 
Here we recall that $\mathbf{\Sigma}(x)=G(x) G(x)^T$ ($x\in \mathbb{R}^p$), and that $\mu_y$ is the disintegration of the ergodic statistical equilibrium $\mu$ ``above'' $y \in V$; see~\cite[Sec.~3]{Chekroun2017a} for more details.
When $h$ is the projection onto $V$, the probability measure $\mu_y$ is the conditional probability of the unobserved variables, contingent upon the value of the observed variable to be $y$. 

Thus if the semigroup property is satisfied for the limiting family of Markov operator $(\mathfrak{T}_\tau)_{\tau\geq 0}$, the reduced RP resonances possess a dynamical interpretation as they relate, through $\overline{\mathcal{A}}$, 
 to the original drift and diffusion terms, $F$ and $G$, respectively. Furthermore, the spectral theory of semigroups applies to $(\mathfrak{T}_\tau)_{\tau \geq 0}$ in this case and one can decompose using~\cite[Theorem 2.5 and Corollary 2.1]{Chekroun2017a}  correlation functions of the reduced state space $V$ (not correlations of the full state space) by using the spectral elements of $\mathfrak{T}_\tau$.
Such correlations are for instance time-lagged cross-correlations between two components of $Y_n$, if dim$(V) \geq 2$.

As a consequence,  still denoting by $\lambda_j$ the eigenvalues of $\overline{\mathcal{A}}$, the correlation function $C_{f,g}(t)$ between two functions $f ,g : V\rightarrow \mathbb{R}$ (i.e.~two observables of the reduced state space $V$) decomposes as

\begin{equation}\label{Eq_corr_decomp}
C_{f,g}(t) =\sum_{j = 1}^\infty e^{\lambda_j t} w_j(f, g)-\langle f\rangle_\mathfrak{m} \langle g\rangle_\mathfrak{m}, 
\end{equation}
%%%%
where $w_j(f,g)=\langle \psi_j^*, f \rangle_\mathfrak{m} \langle g, \psi_j\rangle_\mathfrak{m}$, with $ \psi_j $ denoting the corresponding eigenfunctions of $\overline{\mathcal{A}}$ (in the observed variables only); and $\langle f\rangle_\mathfrak{m}$ denotes the average of $f$ with respect to $\mathfrak{m}$.
Here, $\mathfrak{m}$ is the image by $h$ of the ergodic statistical equilibrium $\mu$ on $V$.
If $h$ consists of the projection onto $V$, then $\mathfrak{m}$ corresponds to the marginal distribution of $\mu$ onto $V$. 

Thus~\eqref{Eq_corr_decomp} shows that in the case $\mathfrak{T}_{\tau+\tau'}=\mathfrak{T}_\tau \mathfrak{T}_{\tau'}$, any correlation function of the reduced dynamics can be exactly decomposed in modes associated with each $\lambda_k$, characterizing thus the mixing properties in the observation space $V$.

To recap, one needs  in practice well chosen reduced variables and that the reduced dynamics encountered is not so intricate as to    
violate $T_{\tau+\tau'}=T_\tau T_{\tau'}$.  Because the checking of the latter condition is not always possible (as it would require a very long temporal dataset),
it is instead very useful to compare  in practice --- based on~\eqref{Eq_corr_decomp} --- correlation functions of the reduced dynamics as estimated using standard correlation estimators, on one hand, and reconstruction formulas involving the spectral elements of $T_\tau$, on the other. 
Indeed if the semigroup property holds, the spectral elements of $T_\tau$ must allow for a decomposition of $C_{f,g}(t)$ as in~\eqref{Eq_corr_decomp} only up to some discretization error related to the mesh size associated with $M$, and the temporal length of the dataset given by $N$.

We describe next how to proceed with such a verification in practice.
In that respect, we first introduce some notations. 
The vector $\bv{m}$ with components given by

\begin{align}
  m_i = \frac{\#\{Y_{n} \in B_i\}}{\#\{Y_{n}
  \in \mathcal{D}\}},\label{eq:estim_reduced_measure}
\end{align}
provides, over the partition $\mathcal{P}_M = \{B_i, \;\, i=1,\cdots, M\}$, the discrete approximation of the density of $\mathfrak{m}$.
It corresponds to the sojourn time density, i.e.~the relative occupancy of a box by the $Y_n$'s, or histogram.
The column vector $\bv{u}$ corresponds to the projection of some function $u$ on the function space generated by the partition $\mathcal{P}_M$, $\bv{u}^\ast$ its conjugate transpose and $\bv{u}'$ simply its transpose.
Furthermore $D(\bv{m}) $ is the diagonal matrix with entries given by the $m_i$'s  ($1\leq i \leq M$), the matrix exponential, $e^{t \bv{\Lambda}_\tau}$, is the diagonal matrix with entries  $e^{t \lambda_k(\tau)}$, and $\bv{\Psi}_\tau $ is the matrix whose columns are given by the eigenvectors of $T_\tau$. 

These notations being clarified, given two observables, $f,g : V \rightarrow \mathbb{R}$, we consider then the function $\mathcal{C}^M_{f, g}(t)$ that corresponds to the discretization of the RHS of~\eqref{Eq_corr_decomp}, namely 

\begin{align}\label{eq:matrixCCFApprox}
  \mathcal{C}^M_{f, g}(t)
  = \left(\bv{f}' D(\bv{m}) \bv{\Psi}_\tau \right)
  e^{t \bv{\Lambda}_\tau}
  \left(\bv{\Psi}^*_\tau{}' D(\bv{m}) \bv{g}\right)
  - (\bv{m}' \bv{f}) (\bv{m}' \bv{g}).
\end{align}
\eqref{eq:matrixCCFApprox} is obtained by discretizing the RHS of~\eqref{Eq_corr_decomp}, in which the inner products involved therein in the weights, $w_j(f,g)$, have been approximated over the partition $\mathcal{P}_M$. Therefore assuming a sufficiently long dataset in time, if $\mathcal{C}^M_{f, g}(t)$ approximates $C_{f,g}(t)$ (as estimated from standard correlation estimators) such that the error shrinks as $M$ increases, one can reasonably infer that $T_{\tau+\tau'}=T_\tau T_{\tau'}$, or at least that this semigroup holds, approximately. This is due to the fact that~\eqref{Eq_corr_decomp} is a necessary condition for $\mathfrak{T}_{\tau+\tau'}=\mathfrak{T}_\tau \mathfrak{T}_{\tau'}$ to hold, for the limiting family of Markov operator $(\mathfrak{T}_\tau)_{\tau\geq 0}$.

The checking of $C_{f,g}(t) \approx \mathcal{C}^M_{f, g}(t)$ does not only nurture a theoretical understanding about the reduced RP resonances as estimated from time series but has important consequences in practice. 
For instance, the theory of~\cite[Sec.~3]{Chekroun2017a} shows that in this case (and if $h$ is a projector), there exists a reduced SDE that approximates (well) the reduced dynamics in $V$, and that this SDE is given by 

\begin{equation}\label{Eq_ReducedSDE}
\d y = \overline{F}(y) \d t + \sigma(y ) \d W_t^V,
\end{equation}
with $\overline{F}(y)=\int F(x) \d \mu_y $, $W_t^V$ is a Brownian motion in $V$, and $\sigma(y)$ satisfies for $1\leq i,j\leq \mbox{dim}(V)$,

\begin{equation}\label{Eq_diffusion_averaged}
(\sigma(y)\sigma(y)^T)_{ij}=\sum_{k=1}^q \overline{G_{ik}(x)G_{jk}(x)}, % \; \; \mbox{for some } d\geq 1,
\end{equation}
where the $G_{lm}$'s are the diffusion coefficients of the SDE~\eqref{Eq_SDE}, and $\overline{(\cdot )}$ still denotes the averaging operation with respect to the disintegrated measure $\mu_y$; see~\cite[Theorem 3.2]{Chekroun2017a}. 

In practice $ \overline{F}$ and $\sigma$, as subordinated to the knowledge of $\mu_y$,  are out of reach except in special cases (see~\cite[Sec.~4]{Chekroun2017a}), rendering at first sight difficult the efficient determination of the reduced SDE~\eqref{Eq_ReducedSDE}.
This is where the analysis of the geometric pattern formed by the reduced RP resonances in the complex plane comes into play to bypass this difficulty.
For instance if a geometric pattern is identified that corresponds to a specific SDE, then this specific SDE gives actually the reduced SDE~\eqref{Eq_ReducedSDE} itself, and thus reduced RP resonances can be used as a diagnosis tool of reduced SDEs.
In that respect,~\cite{Tantet2017b} have analyzed in details the RP resonances associated with the normal form of a Hopf bifurcation subject to additive noise, as a {\it paragon} of nonlinear oscillations in presence of noise.  

% For instance, as explained below in the case of the CZ model, reduced RP resonances are found to be organized along parabolas for certain parameter regimes; discrete parabolas that are not only the signature of a nonlinear oscillation in presence of noise but also indicate that the corresponding reduced dynamics can be reasonably well approximated by a simple Hopf normal form subject to additive white noise~\cite{Tantet2017b}; see also~\cite{KCB18} for applications to wind-driven ocean gyres.    

\br
Complimentary to the checking $C_{f,g}(t) \approx \mathcal{C}^M_{f, g}(t)$ (in order to check $T_{\tau+\tau'}=T_\tau T_{\tau'}$), one could also, given an observable $f:V\rightarrow \mathbb{R}$, look at a similar verification for the power spectrum, $S_f$, in terms of the sojourn time density $\bv{m}$ and the spectral elements of $T_\tau$. This operation consists of analyzing to which level of accuracy the following approximate equality holds 

\begin{align}
	S_f(z)	\approx	-\frac{1}{\pi}
						\left(\bv{f}' D(\bv{m}) \bv{\Psi}_\tau \right)
						\Gamma_\tau(z)
						\left(\bv{\Psi}^*_\tau{}' D(\bv{m}) \bv{f} \right)
						- (\bv{m}' \bv{f})^2, \qquad z\in \mathbb{C},
						\label{eq:matrixPSDApprox}
\end{align}
where the matrix $\Gamma_\tau(z)$ is defined as the diagonal matrix with entries 

\begin{equation*}
(\Gamma_\tau(z))_{kk} = \frac{\Re(\lambda_k(\tau))}{((z - \Im(\lambda_k(\tau)))^2 + \Re(\lambda_k(\tau))^2)}, \qquad 1 \le k \le n. 
\end{equation*}
\er

%%%%%%%%%%%%%%%%%%%CZ Model%%%%%%%%%%%%%%%%%
\section{Results}\label{sec:resultsCZ}

We apply the reduction method (Sect.~\ref{Sec_SL_explained0}) to analyse the RP resonances in the stochastic CZ model along the deterministic Hopf bifurcation.
The CZ model and the introduction of a stochastic forcing via the wind stress is first presented.

\subsection{Stochastic Cane-Zebiak model}\label{sec:StochasticCaneZebiakModel}
	
The CZ model is composed of a 1.5-shallow-water ocean component with an embedded mixed layer coupled to a steady-state linear shallow-water atmosphere, both on an equatorial beta plane (for a detailed description, see~\cite{Zebiak1987a}).
The domain, centered at the equator, is rectangular with western and eastern boundaries at the position of Indonesia and Latin America,
respectively, while it is infinite in the meridional direction.
With such boundary conditions, the spatial fields, such as the Sea Surface Temperature (SST) and thermocline depth, are expanded into spectral basis functions,
with 30 Chebychev polynomials in the zonal direction and 31 Hermite functions in the meridional direction,
for a total of 930 variables per field.
For this study, we use the fully-coupled version of the CZ model~\cite{Vaart2000}.
In this version, not only the anomalies but also the mean fields of the ocean and atmosphere are coupled, thus getting rid of spurious stable solutions that would otherwise exist~\cite{Neelin1995}.
On top of an external wind-stress  forcing (described below), 
the atmosphere responds to changes in the SSTs, while the ocean in turns responds to the total wind-stress.
As a consequence, the long-term dynamics of this model are largely determined
by the non-dimensional coupling parameter $\xi$.
For low values of $\xi$ and a standard choice of the other parameters,
only stable stationary solutions exist.
However, for a coupling larger than the critical value $\xi_c \approx 2.87$,
a supercritical Hopf bifurcation occurs at which a periodic orbit emerges with a period of 3--4 years, reminiscent of ENSO.

As we have seen in the introduction, the impact of a stochastic wind forcing, representing fast atmospheric processes such as westerly wind-bursts, on tropical ocean-atmosphere dynamics, is of great interest for our understanding of ENSO.
In particular,~\cite{Roulston2000} have shown that, in the CZ model, a stochastic wind-forcing is able to excite oscillations with a similar period as ENSO before the criticality, i.e for $\xi < \xi_c$, in particular when the noise is red (see also~\cite[Chapter 8]{Dijkstra2013}).
Here, following the methodology of~\cite{Roulston2000},
a stochastic wind-forcing is added to the deterministic wind-forcing in the CZ model.
The wind-stress is composed of a mean field to which is added a Wiener process with spatial patterns respectively given by the mean and first Empirical Orthogonal Function (EOF) of observed pseudo wind-stress anomalies~\cite{Goldenberg1981}.
As in~\cite{Roulston2000}, the seasonal cycle was removed from the observational records as well as the regression of a linear contribution from observed SST anomalies~\cite{Smith1996}.

\begin{figure}[ht!]
  \centering
  \begin{subfigure}{0.48\linewidth}
    \includegraphics[width=\linewidth]{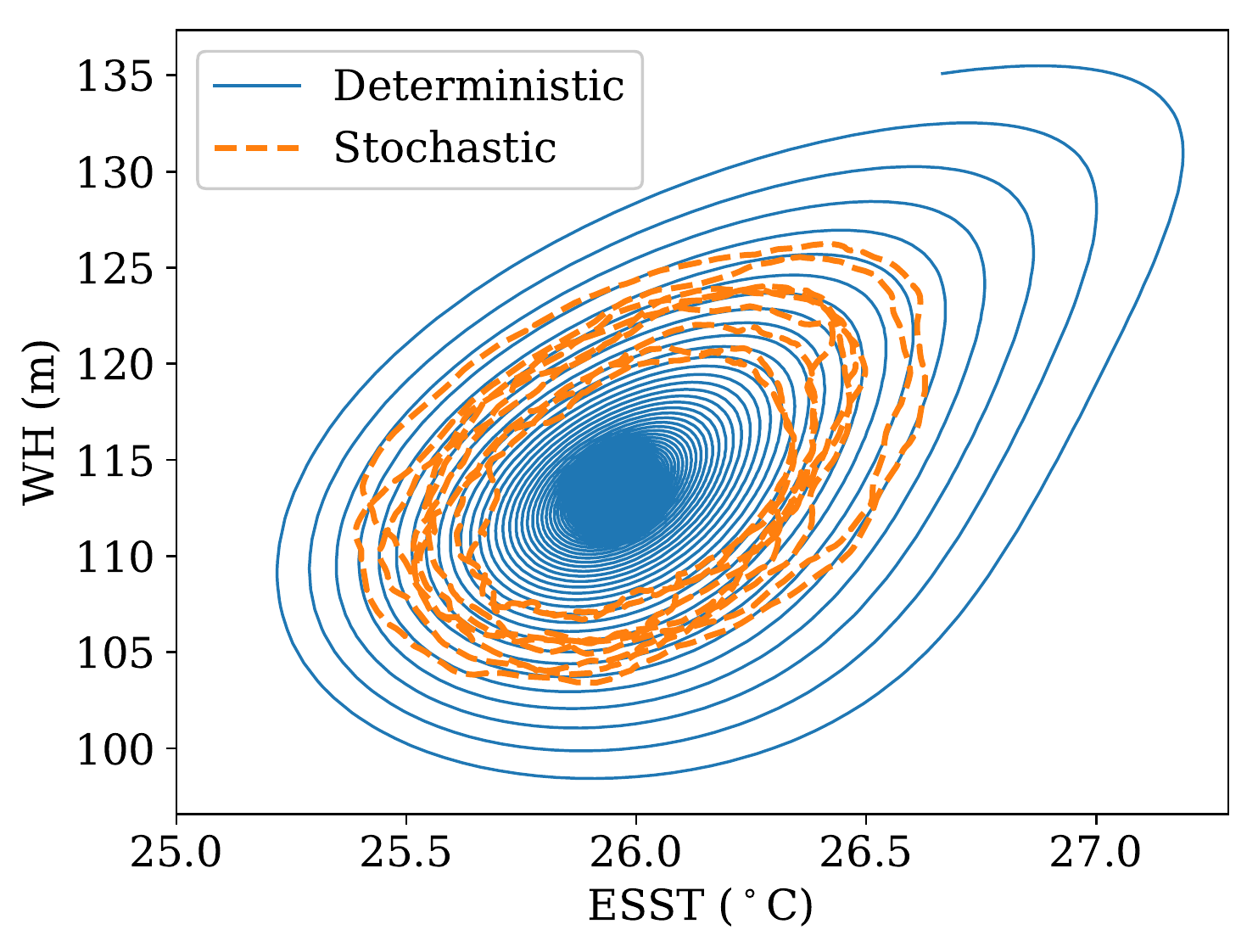}
    \caption{$\xi = 2.85$}
  \end{subfigure}
  \begin{subfigure}{0.48\linewidth}
    \includegraphics[width=\linewidth]{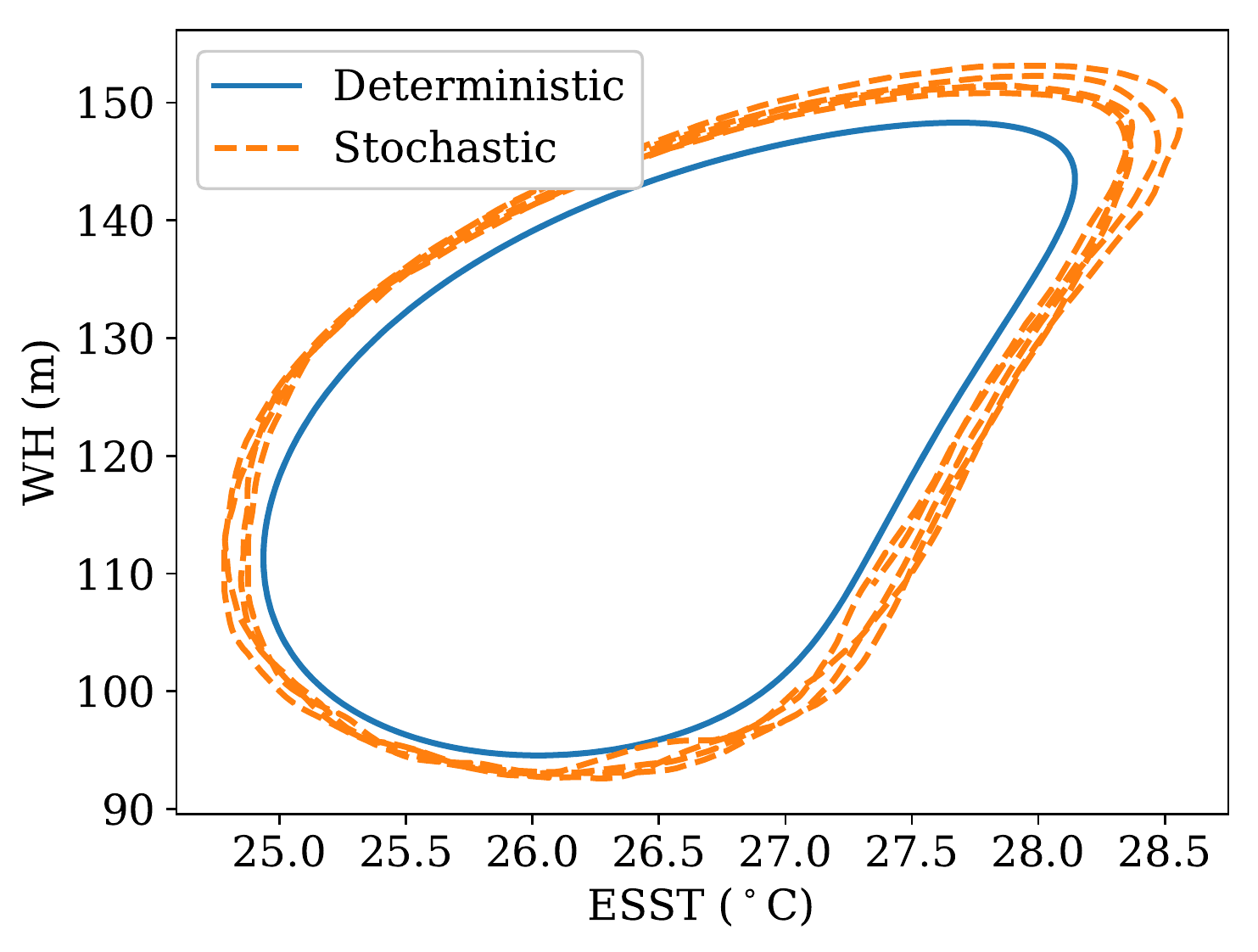}
    \caption{$\xi = 2.90$}
  \end{subfigure}
  \caption{Trajectories projected in the (ESST, WH) plane of CZ model integrations without (plain blue line) and with (dashed orange line) stochastic forcing for $\xi = 2.85$ (left, right before the bifurcation) and $\xi = 2.90$ (right, right after the bifurcation). }\label{fig:time_series}
\end{figure}

We have chosen to model the noise by a Wiener process, i.e.~with equal variance at all frequencies.
The level of the noise, calculated as the spatial average of the standard deviation of the non-dimensional stochastic wind-forcing, was given a value of 0.01.
Model simulations of 6000 years were run with an integration time-step of 5 days,
for different values of the control parameter $\xi$.
In this configuration but without noise, a Hopf bifurcation is observed for $\xi \approx 2.87$.
Deterministic and stochastic trajectories from the CZ model simulations projected in the plane composed of time series of SST to the east of the basin (ESST) and of the thermocline depth to the west of the basin (WH), both at the equator, are represented in Figure~\ref{fig:time_series}.
For a coupling $\xi$ of 2.85 (left panel), the deterministic trajectory (plain blue line) converges to a stationary point, while the stochastic trajectory wanders erratically about this point.
For a coupling of 2.90 (right panel), the deterministic trajectory has converged to the limit cycle, while the stochastic trajectory wanders about this cycle.

\subsection{RP resonances in the CZ model}\label{Sec_mixing_spec_CZ}

Although the CZ model can be considered as a model
of intermediate complexity compared to the state of the art,
it has nonetheless thousands of degrees of freedom.
To analyze the mixing dynamics in state space for the CZ model, we thus apply the reduction method of Section~\ref{Sec_SL_explained0}.
We define an observation operator $h$ on a low-dimensional space $V$.
Based on the knowledge of the physics of the instability
involved in the periodic dynamics (see e.g.~\cite{Vaart2000}),
we choose a two-dimensional space, with
as first component $h_1$ the east-equatorial SST (ESST)
and as second component $h_2$ the west-equatorial thermocline depth (WH).
This plane is the same as the one of Figure~\ref{fig:time_series}. 
This so-defined reduced state space
is then discretized into 100-by-100 grid-boxes
spanning a rectangle of $[-4.5, 4.5]\times[-4.5, 4.5]$ standard deviations.
The transition matrices $T_\tau$ for a transition lag $\tau$ of \SI{4}{years} are estimated from $6000$ year long simulations of the model including the stochastic wind-forcing (with an integration time step of 5 days and a spinup of 100 years removed) for several values of the coupling $\xi$.

The sampling of the time series, the domain, the number of grid-boxes and the lag $\tau$ are key parameters for a proper estimation of the RP resonances.
Here, they are chosen based on robustness tests and so as for the spectral reconstructions of the correlations functions and power spectral densities represented in Figure~\ref{fig:spectrumCZ} to closely fit to the sample estimates of these objects.
In particular, the lag $\tau$ is chosen as a trade-off between estimating the position of the leading resonances well and capturing resonances further away from the imaginary axis.
The robustness to the sampling was tested by comparing eigenvalues estimated from 10 transition matrices estimated leaving one tenth of the original time series out.
The 10 spectra being virtually identical (not shown here), the results can be considered robust to the sampling.
Moreover, the sampling frequency of the time series of 5 days was kept, while a temporal smoothing may be needed in applications for which the physical phenomenon of interest needs to be isolated from irrelevant fast processes.
The convergence of the eigenvalues estimates for an increasing number of grid boxes per dimension (from 10 to 200 with a step of 10) and for an increasing lag (from 1~month to 100~months with a step of 1 month) was also tested.
The number of grid boxes and the lag were chosen so as for the relative change in the real parts of the first 10 eigenvalues from one step to the other to be less than $1\%$.
For more details on the estimation procedure, see~\cite{Chekroun2014},~\cite{tantet_early_2015},~\cite{tantet_crisis_2018} and~\cite{Chekroun2017a}.

\begin{figure}
	\centering
        \begin{subfigure}[b]{0.48\textwidth}
		\includegraphics[width=\textwidth]{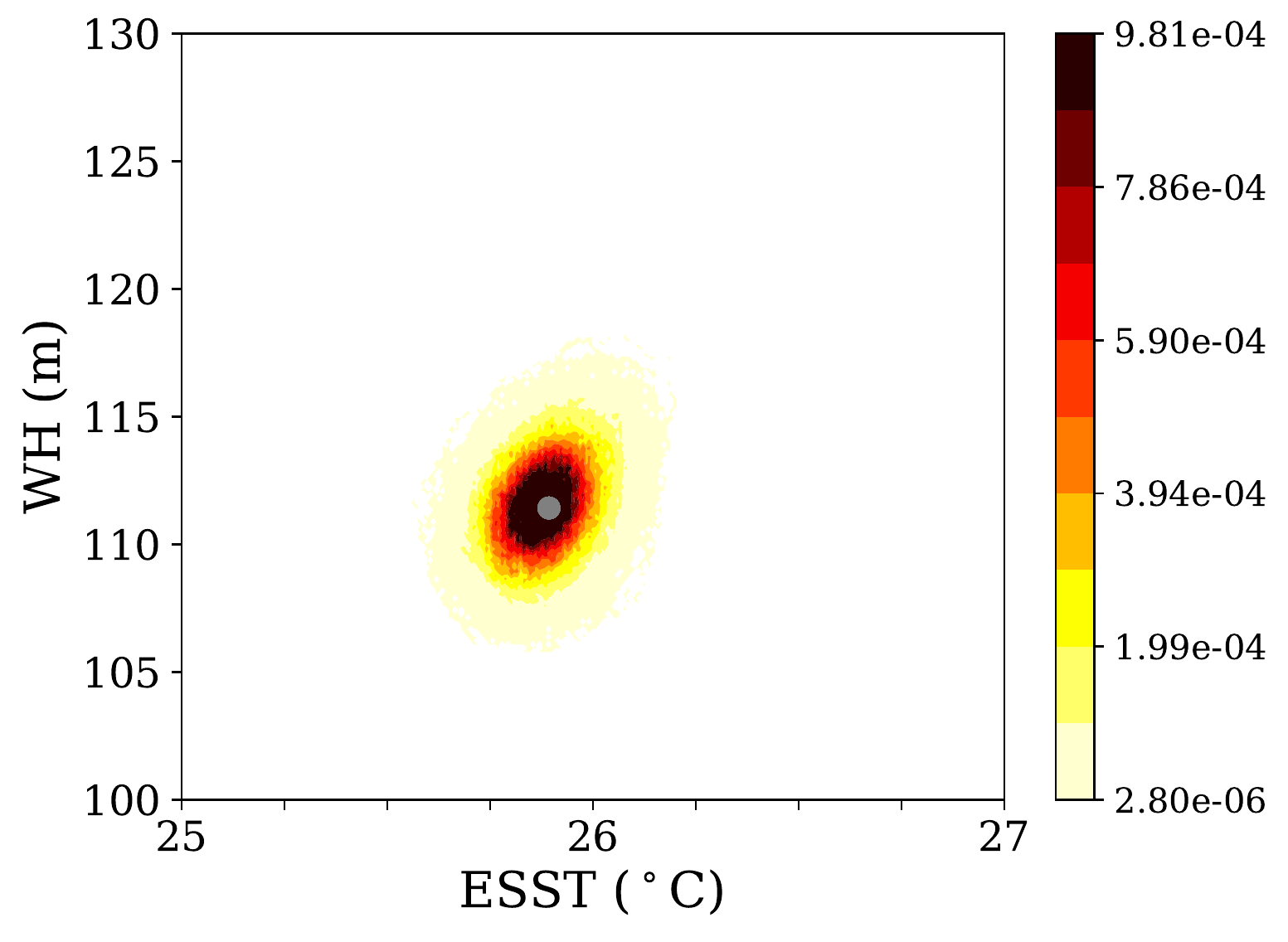}
		\caption{}
	\end{subfigure}
        \begin{subfigure}[b]{0.48\textwidth}
		\includegraphics[width=\textwidth]{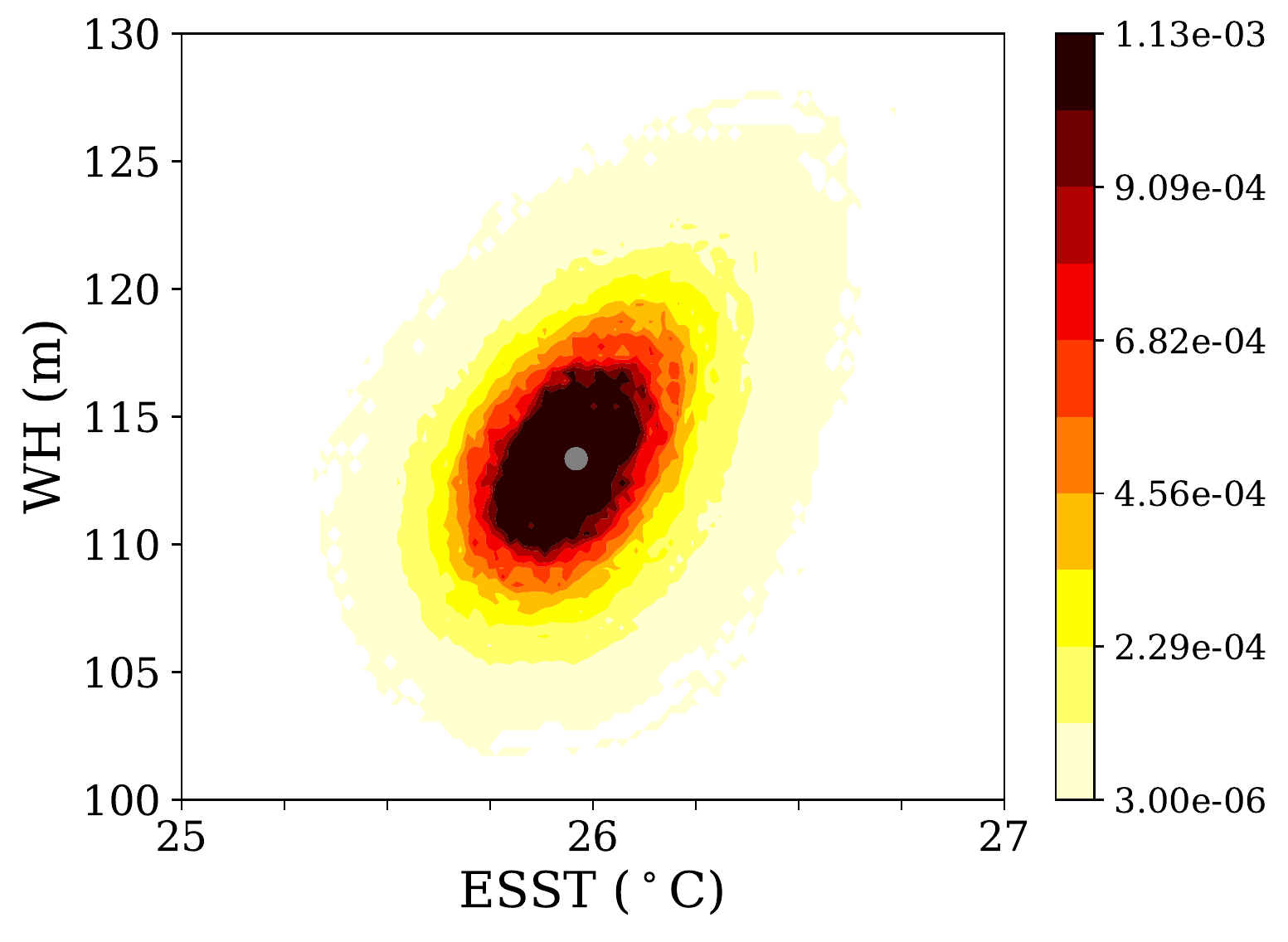}
		\caption{}
	\end{subfigure}\\
        \begin{subfigure}[b]{0.48\textwidth}
		\includegraphics[width=\textwidth]{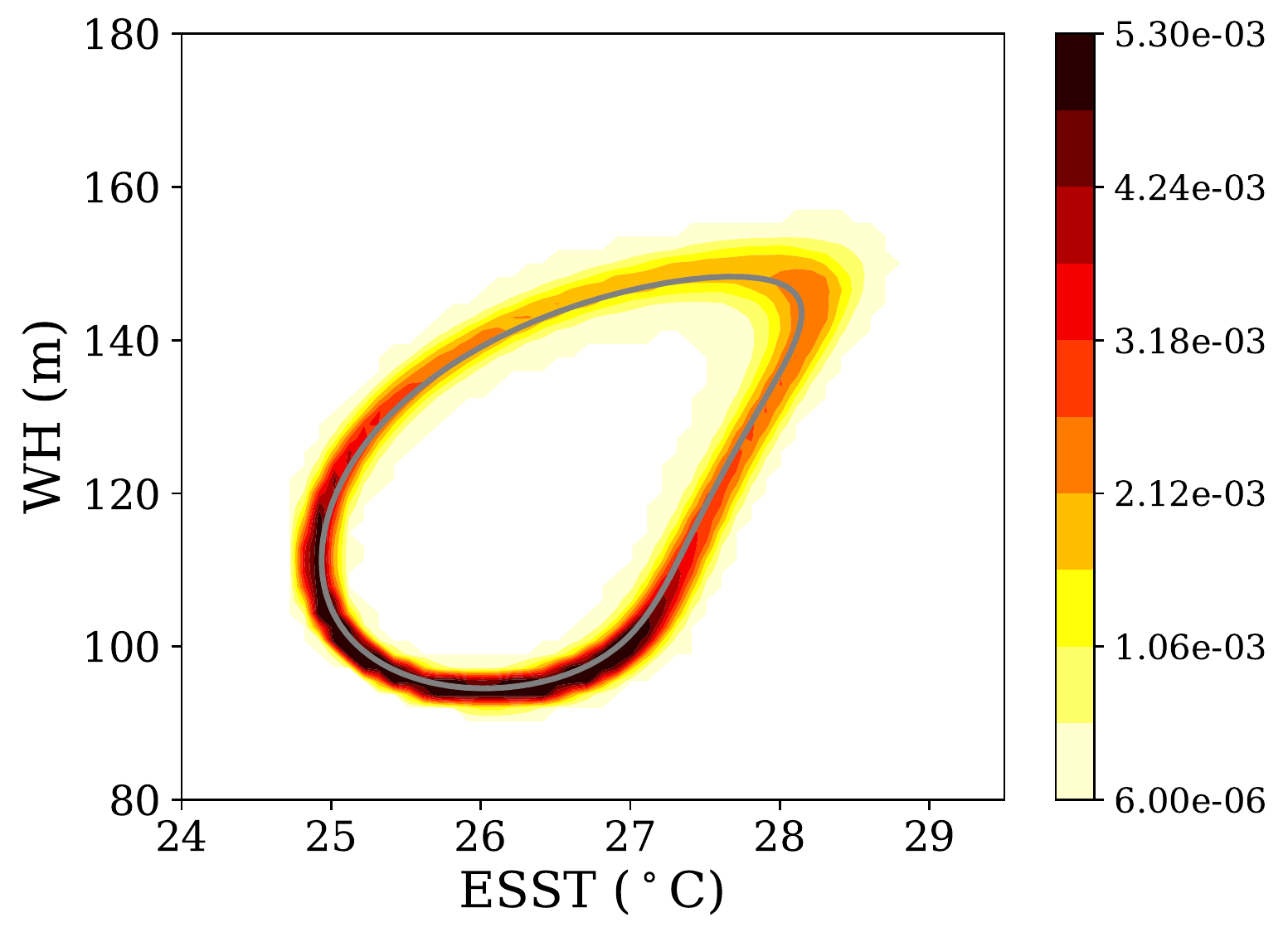}
		\caption{}
	\end{subfigure}
        \begin{subfigure}[b]{0.48\textwidth}
		\includegraphics[width=\textwidth]{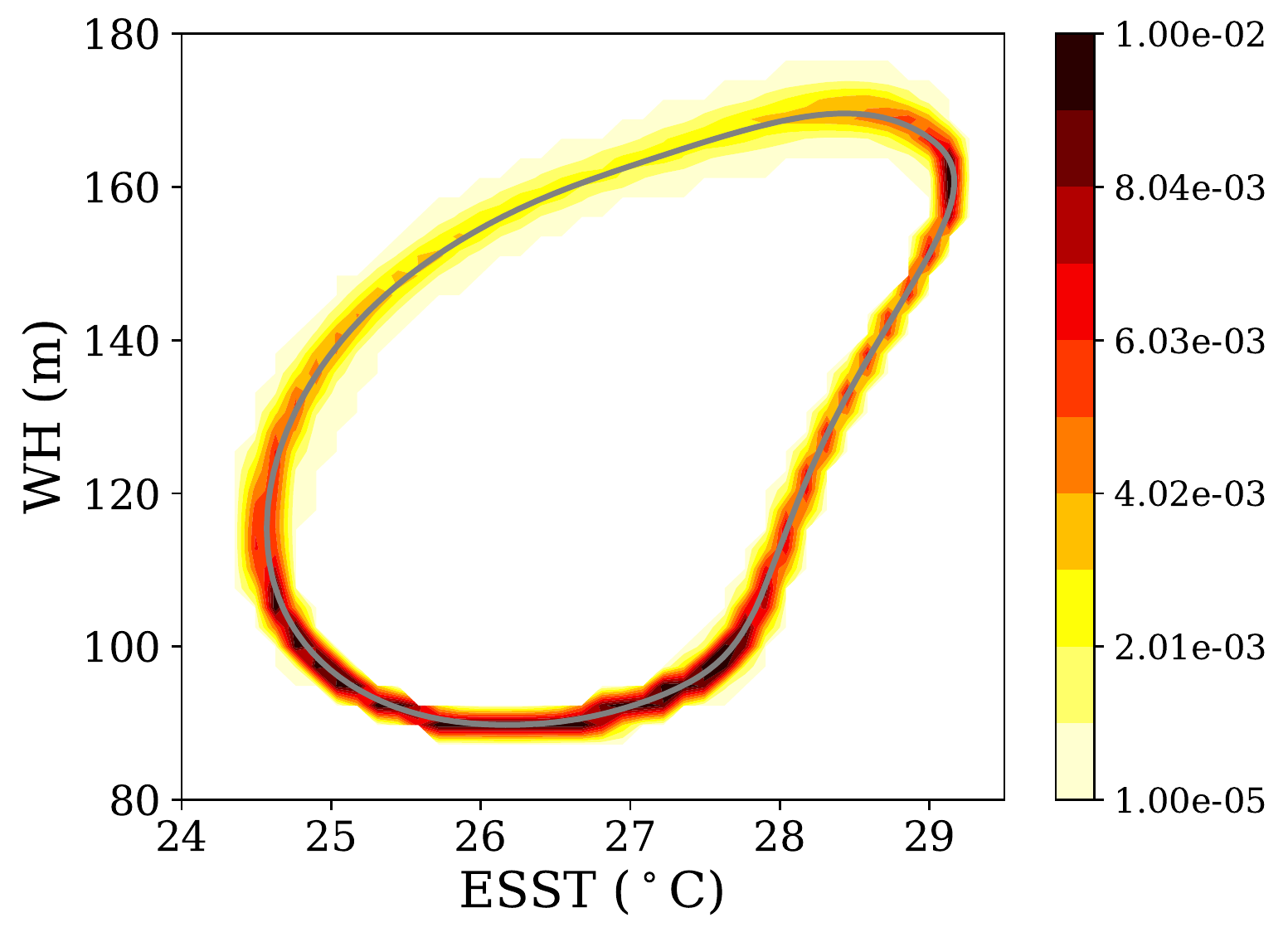}
		\caption{}
	\end{subfigure}
	\caption{Histogram (shaded contours) in the (ESST, WH) plane (same as Fig.~\ref{fig:time_series}) for the stochastic CZ model with a coupling $\xi$ of
          (a)  2.80 (before the deterministic Hopf bifurcation),
          (b) 2.85 (right before the bifurcation),
          (c) 2.90 (right after the bifurcation),
          (d) 2.95 (after the bifurcation).
          The stable stationary point and the limit cycle from the deterministic version of the model and the same values of the coupling are also represented in grey.
          The histogram coincides with the leading left eigenvector of the transition matrix $T_\tau$.
          The scale changes between panels (a,b) and (c,d).}\label{fig:histCZ}
\end{figure}

\begin{figure}
	\centering
        \begin{subfigure}[b]{0.35\textwidth}
          \includegraphics[width=\textwidth]{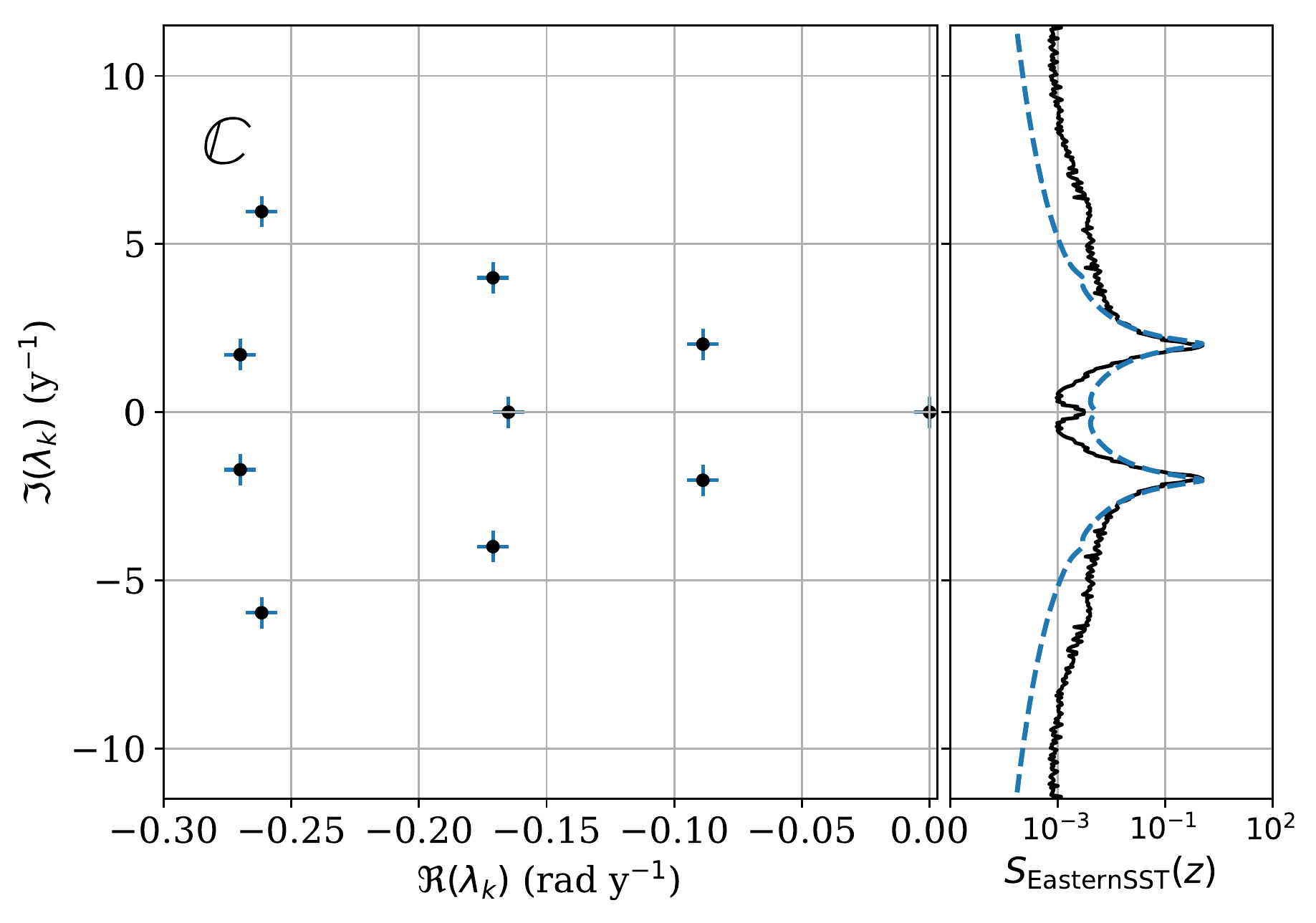}
          \caption{}
	\end{subfigure}
        \begin{subfigure}[b]{0.35\textwidth}
          \includegraphics[width=\textwidth]{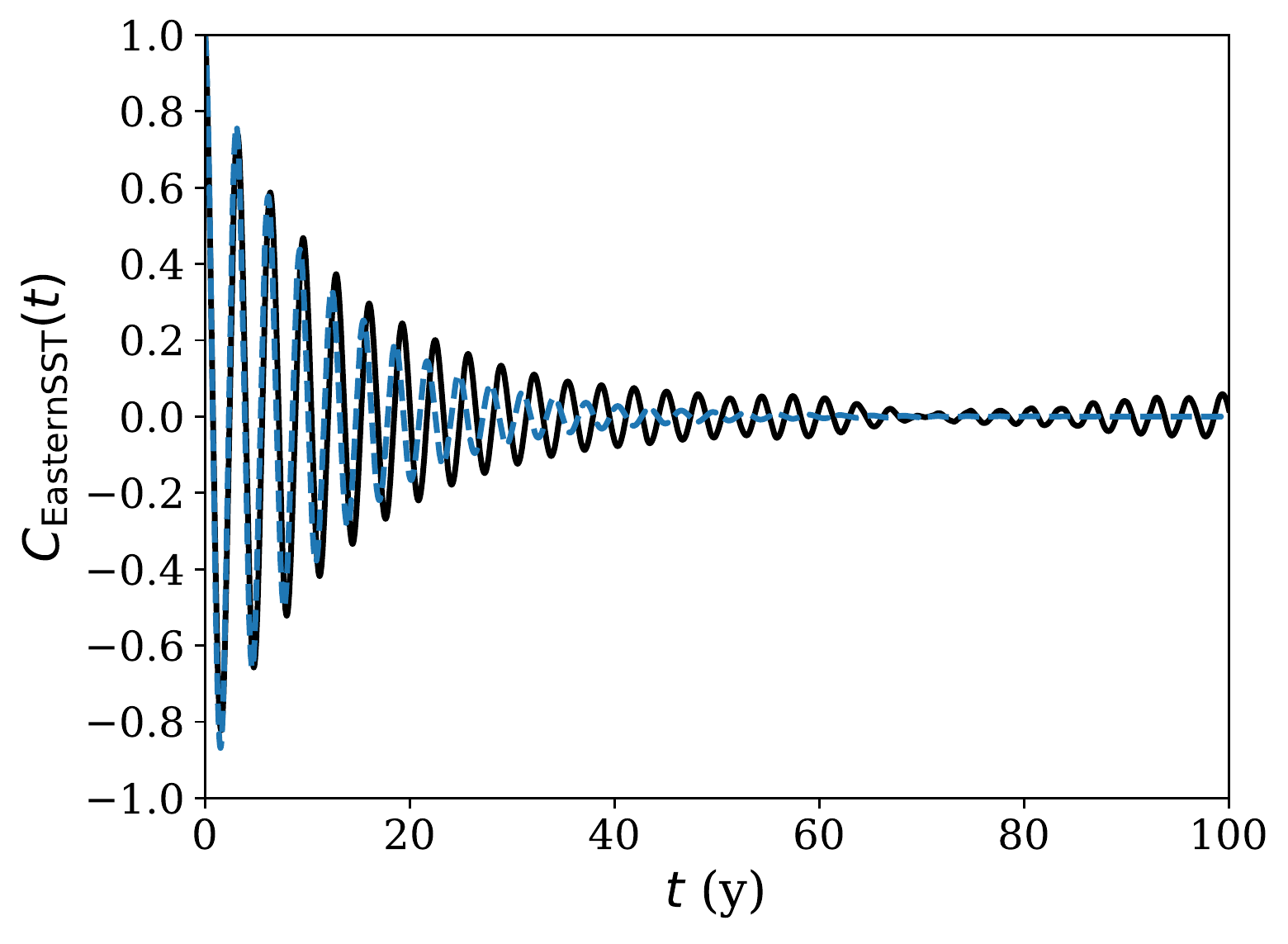}
          \caption{}
	\end{subfigure}\\
        \begin{subfigure}[b]{0.35\textwidth}
          \includegraphics[width=\textwidth]{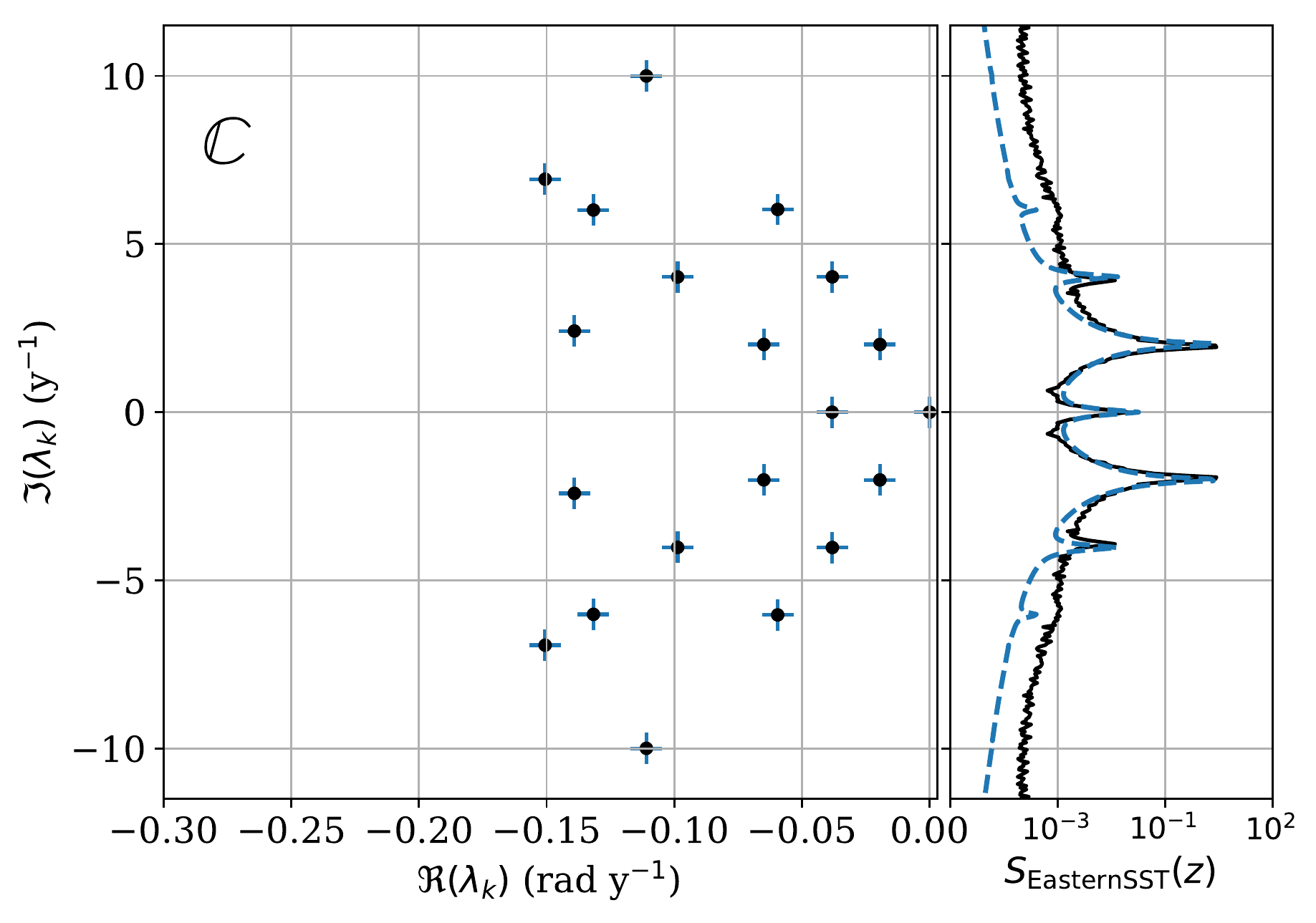}
          \caption{}
	\end{subfigure}
        \begin{subfigure}[b]{0.35\textwidth}
          \includegraphics[width=\textwidth]{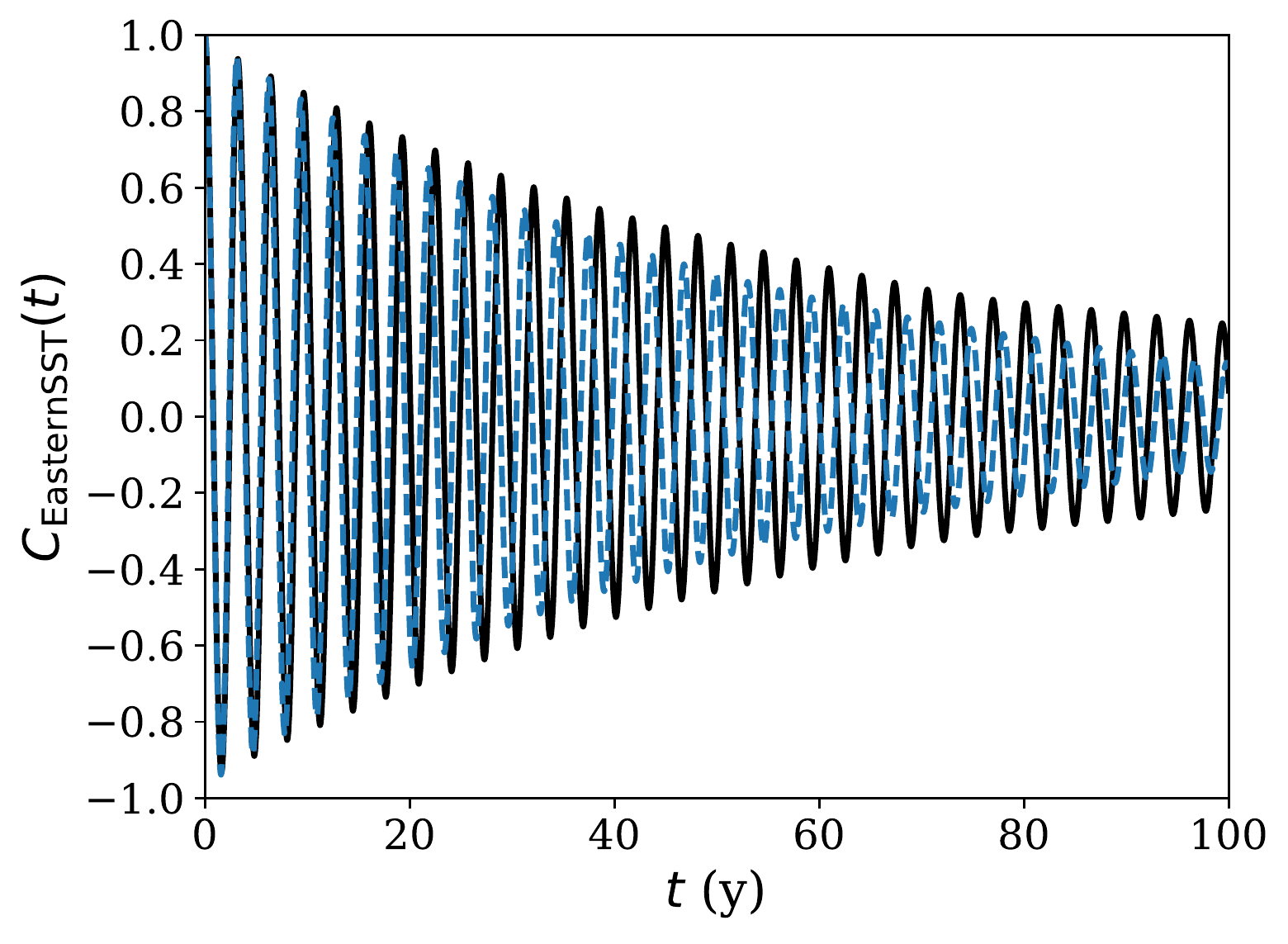}
          \caption{}
	\end{subfigure}\\
        \begin{subfigure}[b]{0.35\textwidth}
          \includegraphics[width=\textwidth]{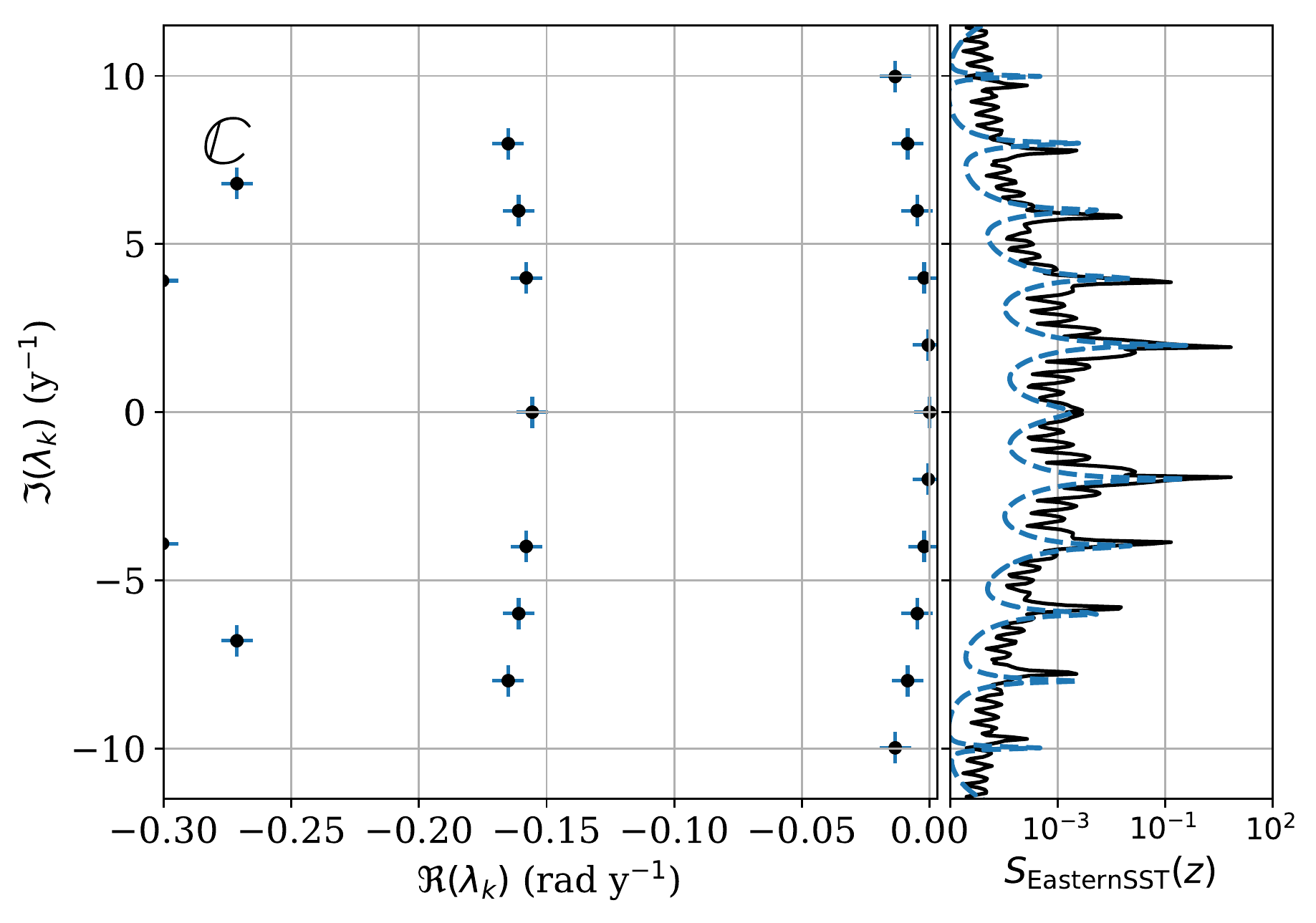}
          \caption{}
	\end{subfigure}
        \begin{subfigure}[b]{0.35\textwidth}
          \includegraphics[width=\textwidth]{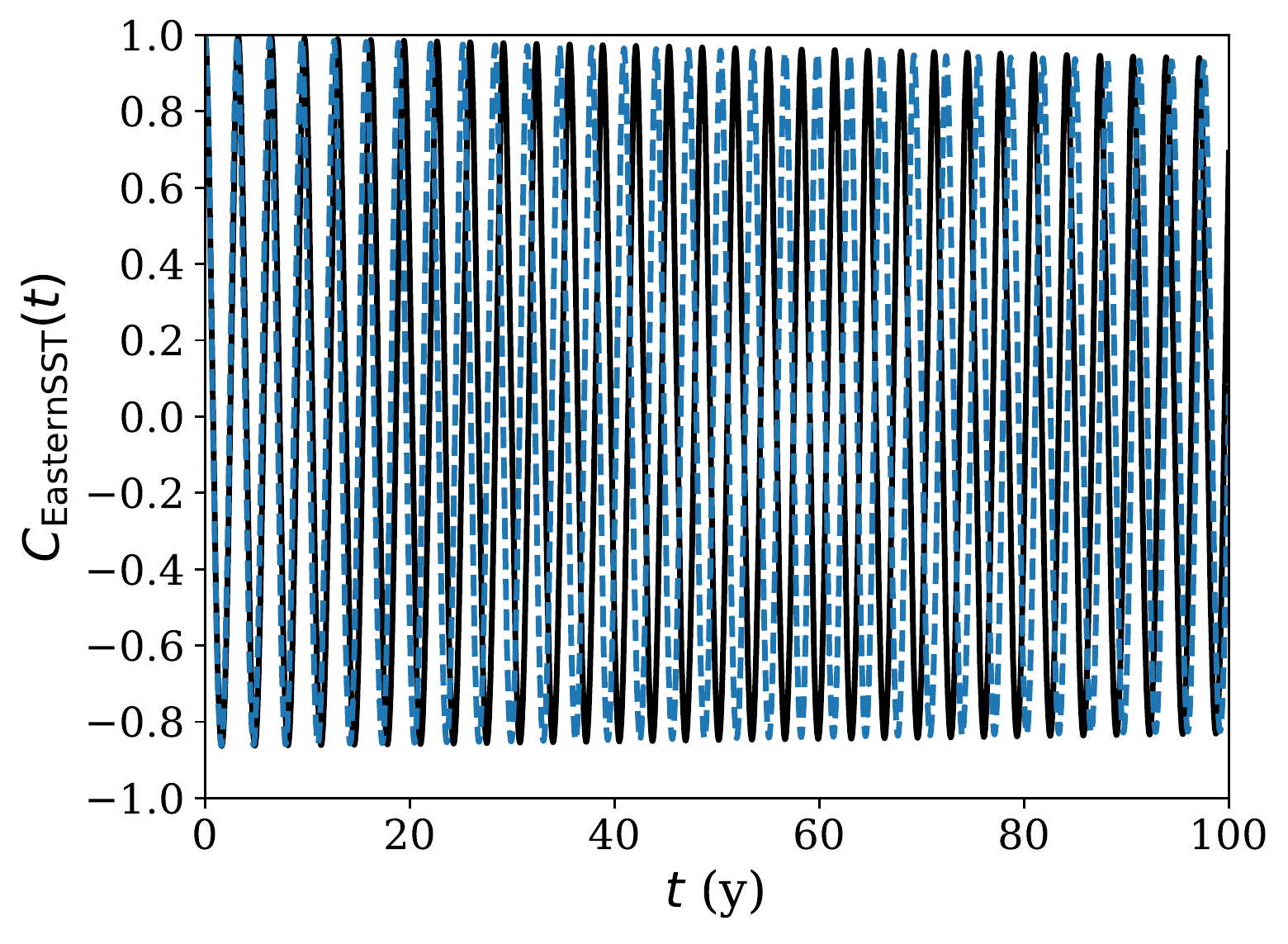}
          \caption{}
	\end{subfigure}\\
        \begin{subfigure}[b]{0.35\textwidth}
          \includegraphics[width=\textwidth]{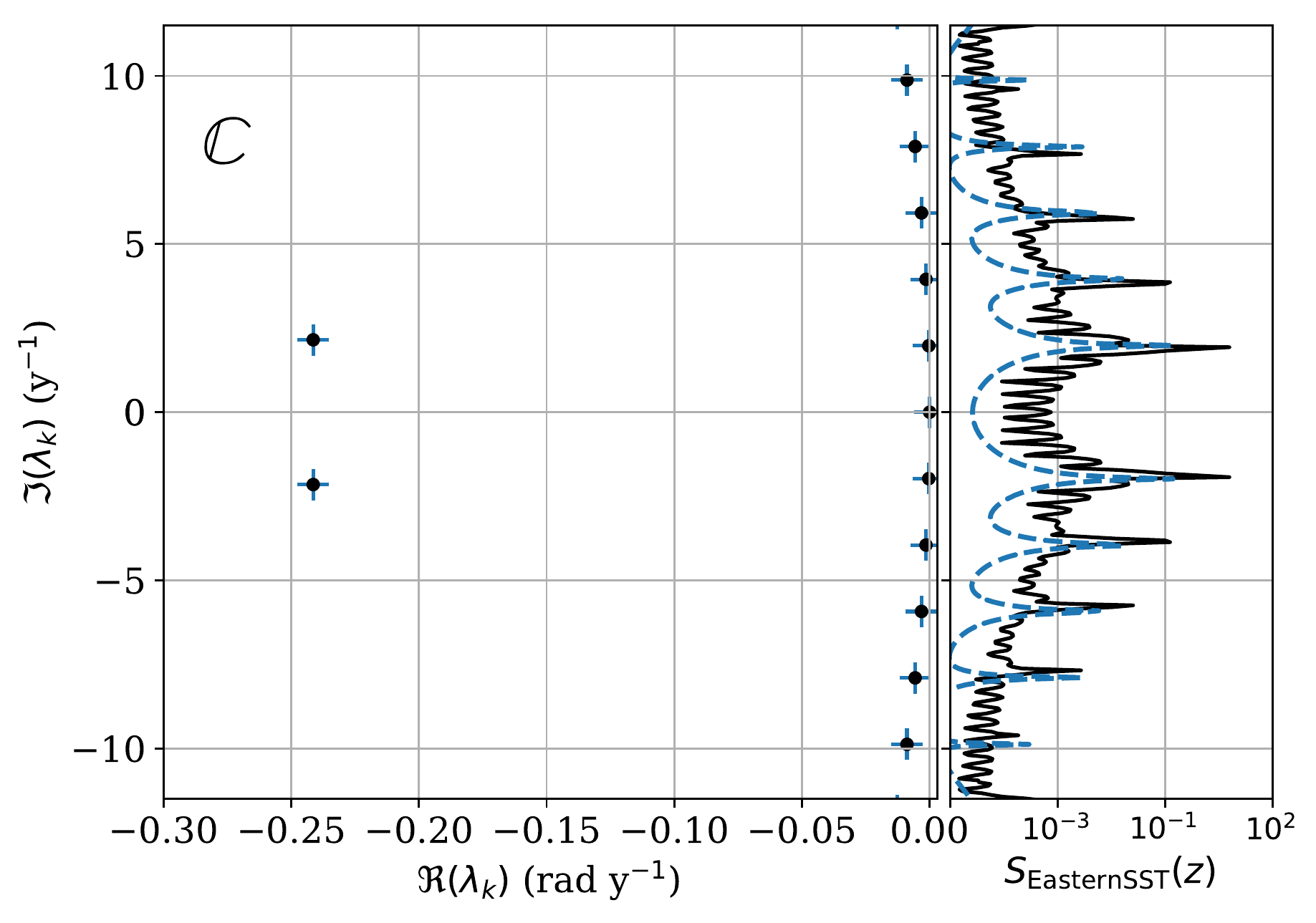}
          \caption{}
	\end{subfigure}
        \begin{subfigure}[b]{0.35\textwidth}
          \includegraphics[width=\textwidth]{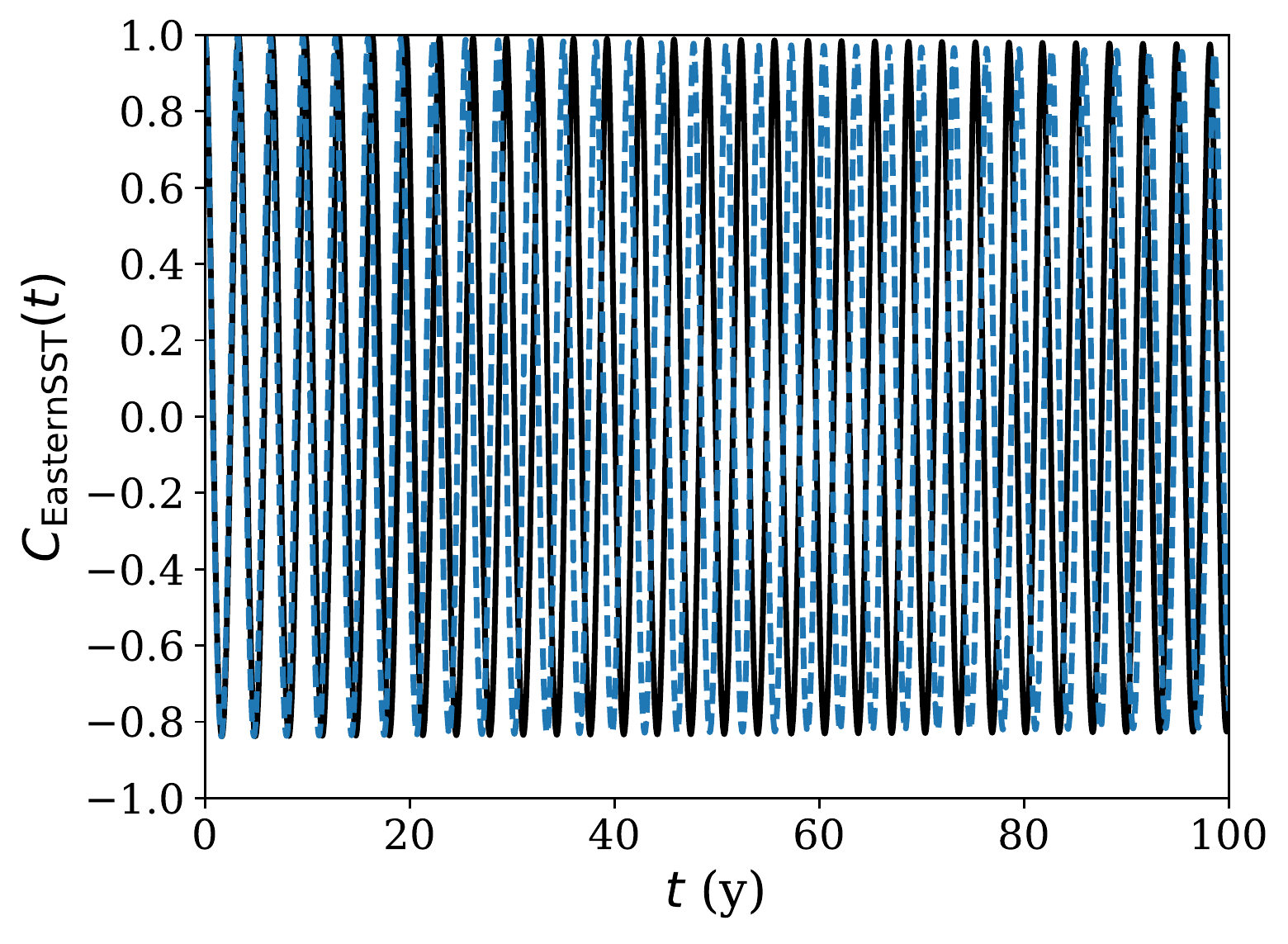}
          \caption{}
	\end{subfigure}\\
	\caption{Leading reduced RP resonances (left part of left panels),
	power spectral density $\psd$ (right part of left panels)
	and correlation function $\acf$ (right panels)
	for the stochastic CZ model with a coupling $\xi$ of
	(a) 2.80 (before the deterministic Hopf bifurcation),
	(b) 2.85 (right before the bifurcation),
	(c) 2.90 (right after the bifurcation),
	(d) 2.95 (after the bifurcation).
	The position of the eigenvalues $\lambda_i$ is marked by a plus with an overlapping disc indicating their relative weight $w_i$, in logarithmic scale.
        Eigenvalue estimates with an eigen-condition number $\kappa_i = |\langle \psi_j^*, \psi_j\rangle_\mathfrak{m}| / (\|\psi_j^*\|_\mathfrak{m}\|\psi_j\|_\mathfrak{m})$ exceeding 5 are deemed not robust and are not represented.
        Here, $\psi_k(\tau)$ and $\psi_k^*(\tau)$ are respectively the right and left $k$th-eigenvectors of the transition matrix $T_\tau$ and $\bv{m}$~\eqref{eq:estim_reduced_measure} is the estimate of the marginal measure $\mathfrak{m}$.
	The spectral reconstructions of the power spectral density and of the auto-correlation function of the ESST are represented by a dashed red line, while the black lines give their sample estimates.}
	\label{fig:spectrumCZ}
\end{figure}

Figure~\ref{fig:histCZ} represents histograms the sojourn time density $\bv{m}$ estimated from these time series.
The latter give an estimate of the marginal distribution $\mathfrak{m}$ of $\mu$ on the reduced space.
For small values of the coupling (panels (a) and (b)), the histogram is spread about the stationary point.
This is expected from Figure~\ref{fig:time_series} where realizations fluctuate erratically around the stationary point, never passing twice through the same point, with a large spread away from the stationary point.
Instead, as the coupling is increased (panels (c) and (d)), the system enters an almost-periodic regime so that the histogram is spread about the deterministic limit cycle.

In the left part of the left panels of Figure~\ref{fig:spectrumCZ}, we represent the leading reduced RP resonances calculated from transition matrices estimated following the methodology of Sect.~\ref{Sec_SL_explained0}, for the same values of $\xi$ as in Fig.~\ref{fig:histCZ}.

The choice of a lag $\tau$ as large as 4 years may come as a surprise, as the oscillations in the CZ model are between 3 and 4 years.
Care should indeed be taken when computing the imaginary part of generator eigenvalues (the resonances) from eigenvalues of a transition matrix at a given lag $\tau$ by applying formula~\eqref{Eq_lambda}, since it involves taking the principal part of the argument of the latter.
This would result in clipping imaginary parts of true resonances larger than $\pm \pi / \tau$ (Nyquist angular frequency) and thus prevent estimating the harmonics in the power spectra associated with these resonances.
In order to avoid this aliasing effect, we instead compute the reduced RP resonances from the ratio of the corresponding eigenvalues of transition matrices at lags $\tau$ and $\tau + \Delta \tau$, with $\tau = \SI{4}{years}$ and $\Delta \tau = \SI{0.5}{months}$.
In other words, we estimate the reduced RP resonances as
\begin{equation}\label{Eq_lambda_ratio}
  \lambda_k(\tau)
  = \frac{1}{\Delta \tau} \log \frac{\zeta_k(\tau + \Delta \tau)}{\zeta_k(\tau)}
  = \frac{1}{\Delta \tau} \left(\log \left|\frac{\zeta_k(\tau + \Delta \tau)}{\zeta_k(\tau)}\right|
  +  i\arg \frac{\zeta_k(\tau + \Delta \tau)}{\zeta_k(\tau)}\right),
  \qquad 1\leq k\leq M.
\end{equation}
Doing so, the folding is limited to imaginary parts larger than $\pm \SI{2 \pi}{months^{-1}}$ and the transition lag may be chosen larger, so as for the real parts of the leading eigenvalues to be better estimated.
To pair eigenvalues of the transition matrices $T_\tau$ and $T_{\tau + \Delta \tau}$ corresponding to the same resonance, the distance between the associated eigenvectors is minimized.
Moreover, the eigen-condition numbers $\kappa_k$'s (see caption) are used to assess the robustness of the corresponding eigenvalues to numerical errors and filter out the least robust ones.

In addition, reconstructions of the power spectral density (right part of left panels) and of the auto-correlation function (right panels) of the ESST are represented by a dashed red line.
They are obtained from the leading eigenvalues and eigenvectors by applying the spectral decompositions~\eqref{eq:matrixPSDApprox} and~\eqref{eq:matrixCCFApprox}, respectively.
The black lines instead give the sample estimates~\cite{VonStorch1999b} of these objects.

One can see that, for small values of the coupling $\xi$, panel (a), the reduced eigenvalues are arranged in a triangular structure.
For large values of $\xi$, panel (g), the reduced eigenvalues are instead arranged in a parabolic structure.
For intermediate values of the coupling, panels (c) and (e), the reduced eigenvalues smoothly move from one arrangement to the other.
A first result is that, even for a low value of the coupling parameter, panels (a) and (b), --- for which only a stationary point would exist without noise --- a broad power spectral peak associated with oscillations in the correlation function is visible, reminiscent of that of the observational record of ENSO.
This peak is due to the presence of resonances with non-zero imaginary part close to the imaginary axis, even below the deterministic bifurcation.
This explains the noise-induced oscillations identified by~\cite{Roulston2000}
in a similar configuration of the CZ model.

Importantly, as $\xi$ is increased, the spectral gap between the eigenvalue zero
and the leading secondary eigenvalues decreases.
This shrinkage of the spectral gap is responsible for the slowing down
of the decay of correlations (shown on the right panels, see below),
as well as for the sharpening of the power-spectral peaks (right part of left panels).
The position of these peaks coincide with the harmonics of the periodic orbit.
As $\xi$ is further increased, there is both a reduction of the gap between the imaginary axis and the first line of eigenvalues (i.e.~the one passing through 0) and an increase of the gap between the latter and the second line of eigenvalues,
as visible from the  transition from panel (c) to panel (e) (the second line of eigenvalues is not visible in panel (g) since only the leading eigenvalues are represented in order to filter out the high harmonics from the first line which are folded).
These phenomena may respectively be interpreted from the fact that, with the increase of the stability of the limit cycle, diffusion along the limit cycle weakens and contraction towards the limit cycle strengthens.

Moreover the spectral reconstructions of the power spectral density (right part of left panel, dashed red line) and the auto-correlation function (right panel, dashed red line) of the ESST from a single transition matrix at the chosen lag $\tau$ match relatively well with the corresponding sample estimates (thin black line).
This suggests that the slow dynamics in the reduced state space considered here is little impacted by memory effects induced by the reduction and that the semigroup property of the Markov semigroup is relatively well preserved; see~\cite{chen2016diversity} for a discussion of the role memory effects in reduced state space for ENSO modeling. 

\begin{figure}
	\centering
	\begin{subfigure}[b]{0.32\textwidth}
		\includegraphics[width=\textwidth]{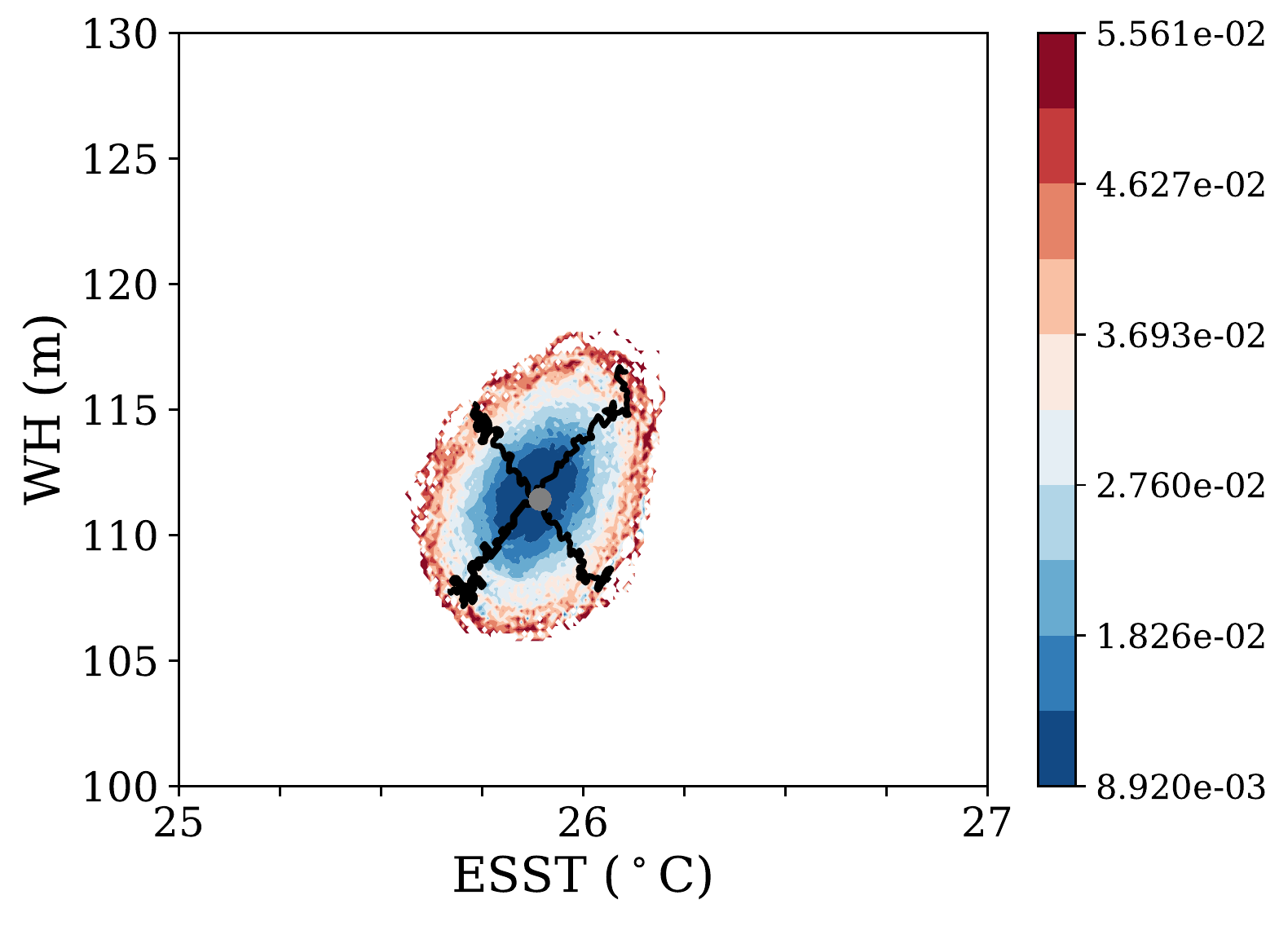}
		\caption{}
	\end{subfigure}
	\begin{subfigure}[b]{0.32\textwidth}
		\includegraphics[width=\textwidth]{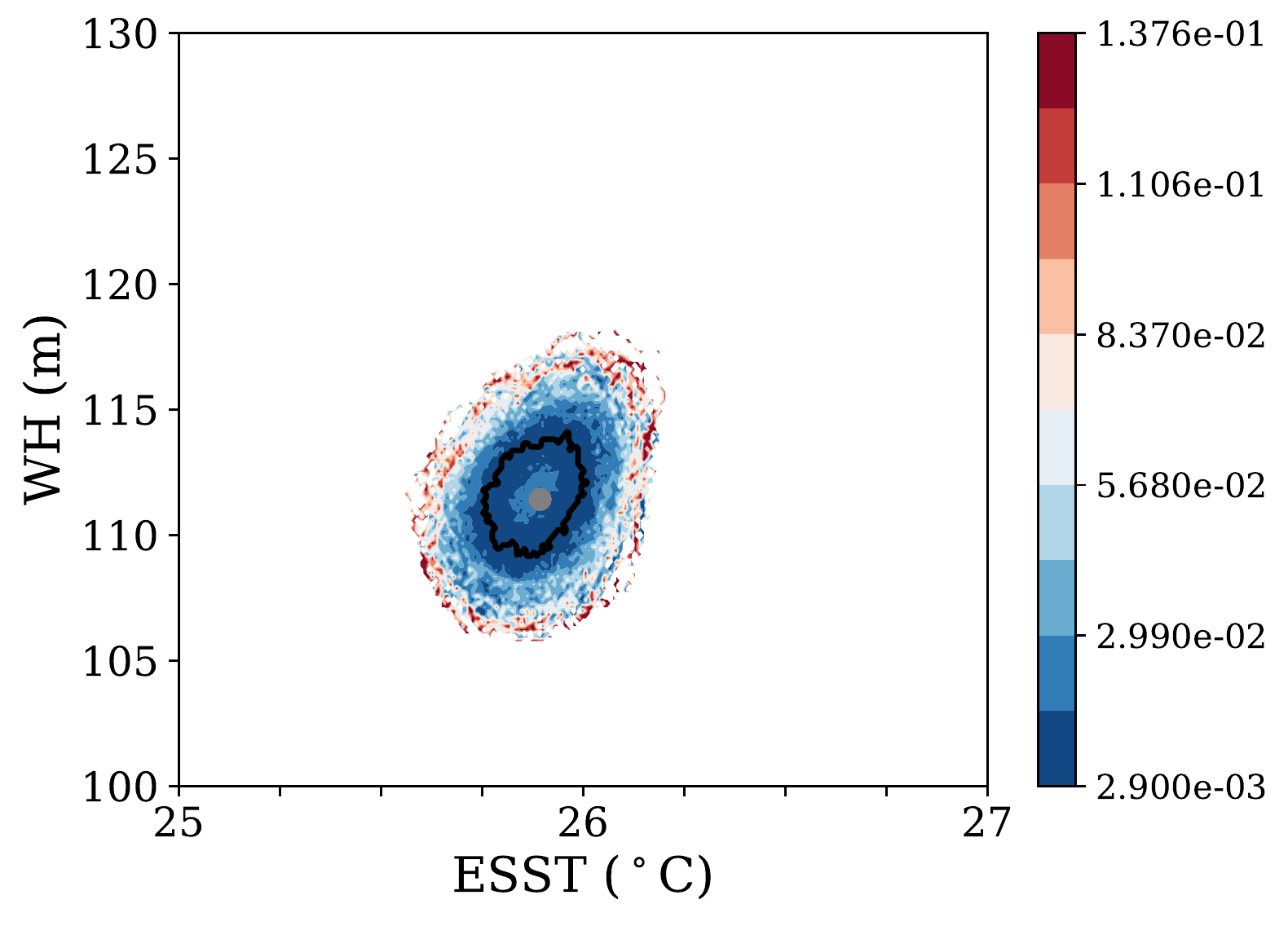}
		\caption{}
	\end{subfigure}
	\begin{subfigure}[b]{0.32\textwidth}
		\includegraphics[width=\textwidth]{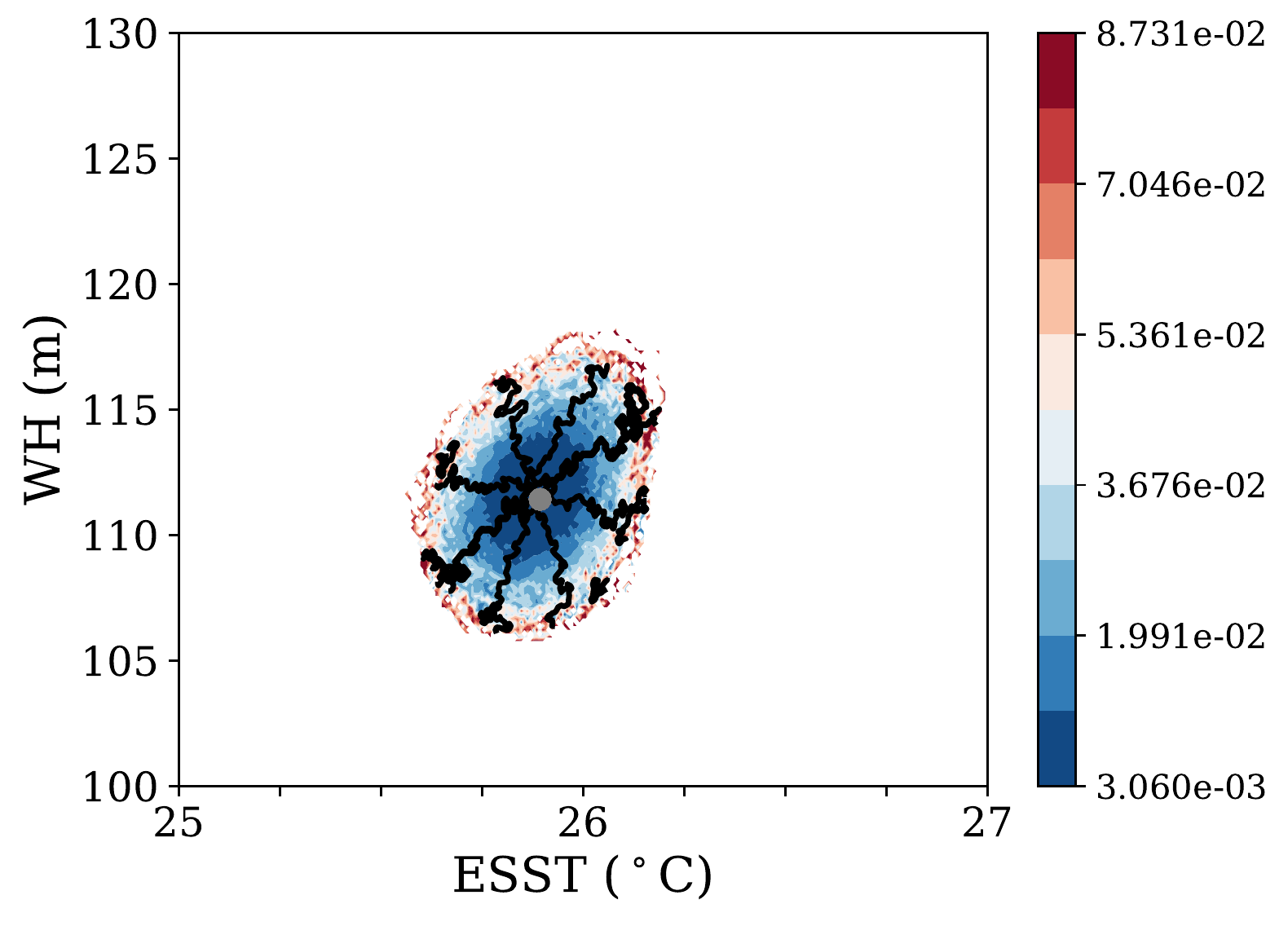}
		\caption{}
	\end{subfigure}\\
	\begin{subfigure}[b]{0.32\textwidth}
		\includegraphics[width=\textwidth]{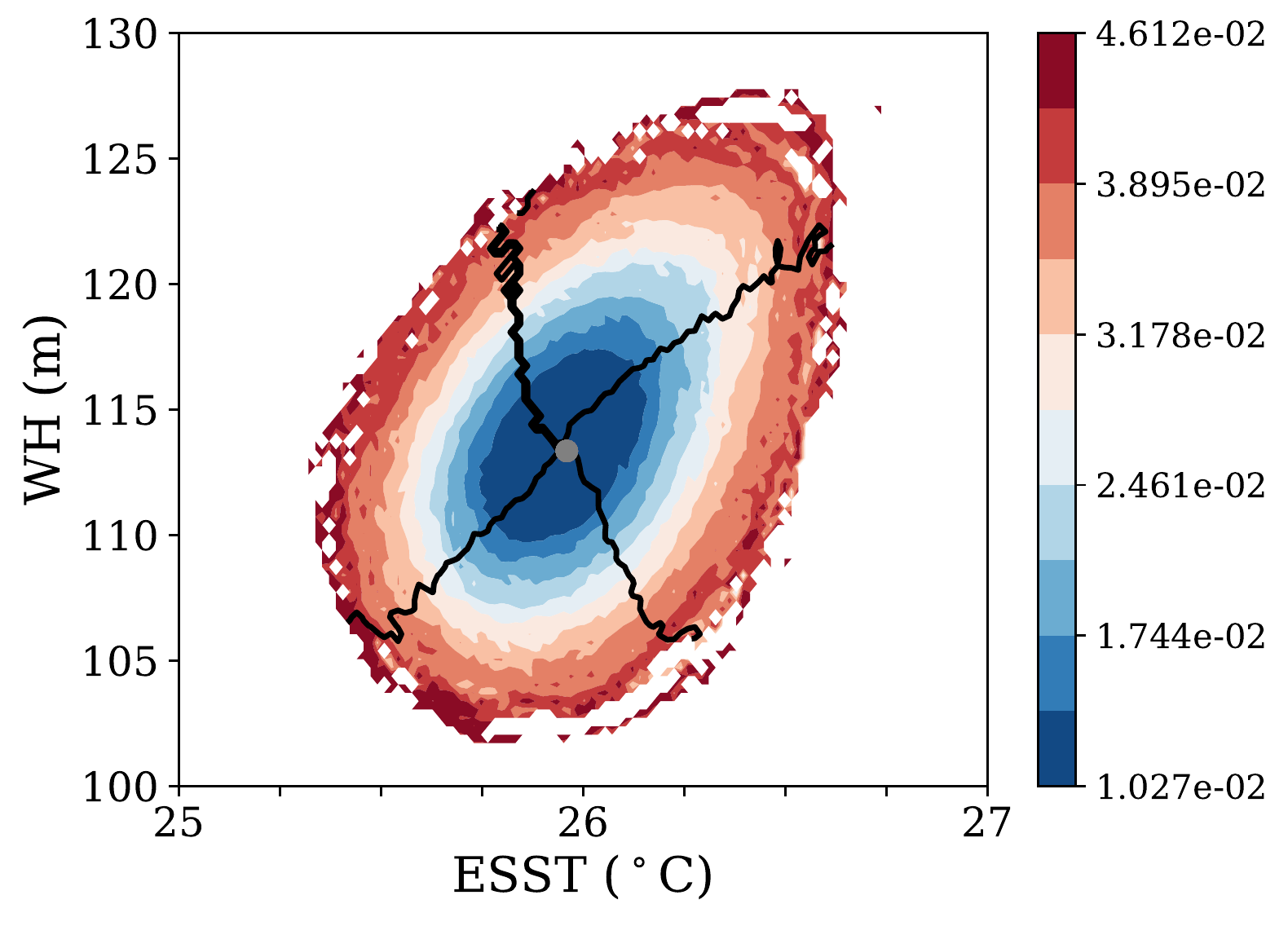}
		\caption{}
	\end{subfigure}
	\begin{subfigure}[b]{0.32\textwidth}
		\includegraphics[width=\textwidth]{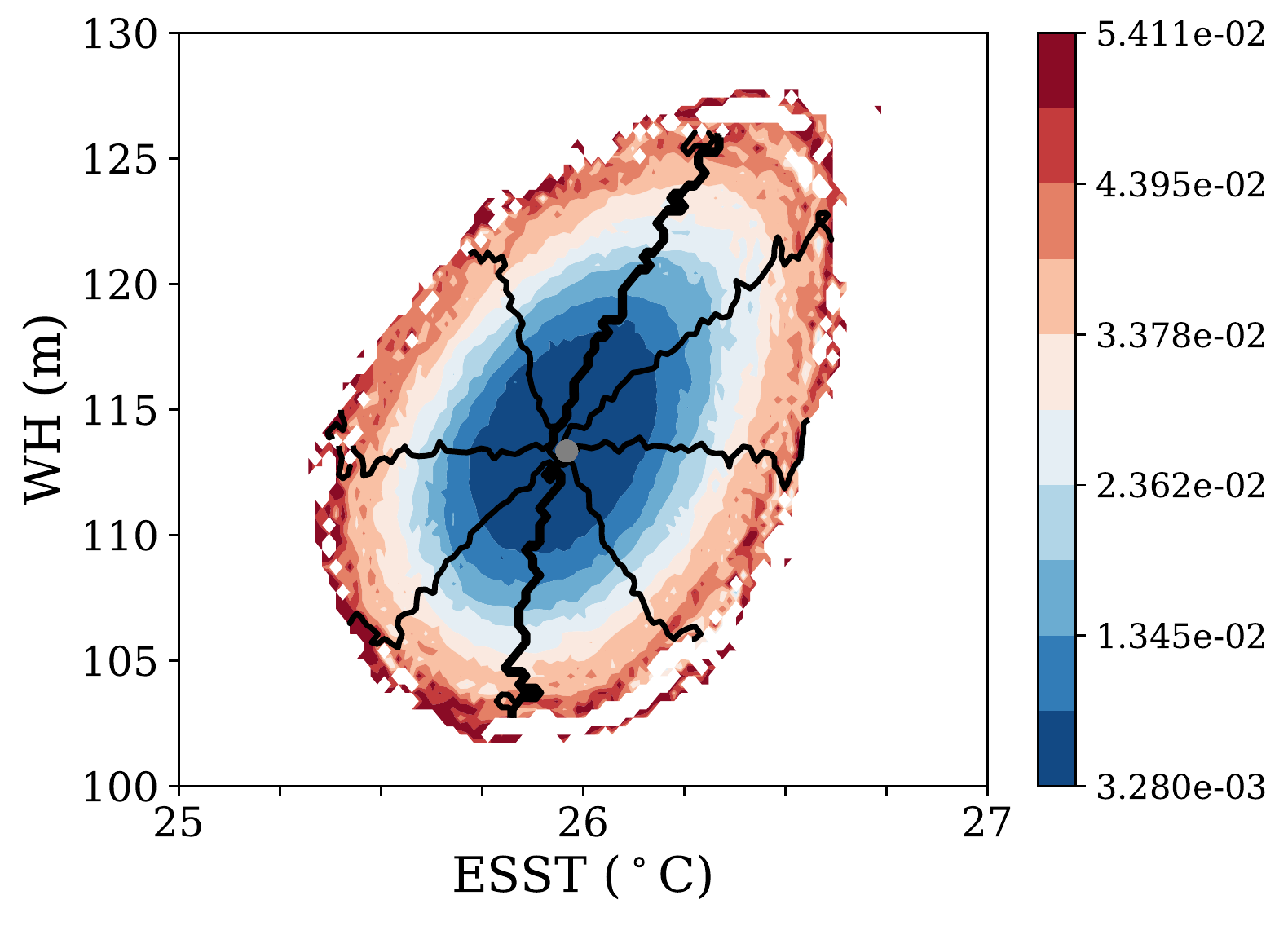}
		\caption{}
	\end{subfigure}
	\begin{subfigure}[b]{0.32\textwidth}
		\includegraphics[width=\textwidth]{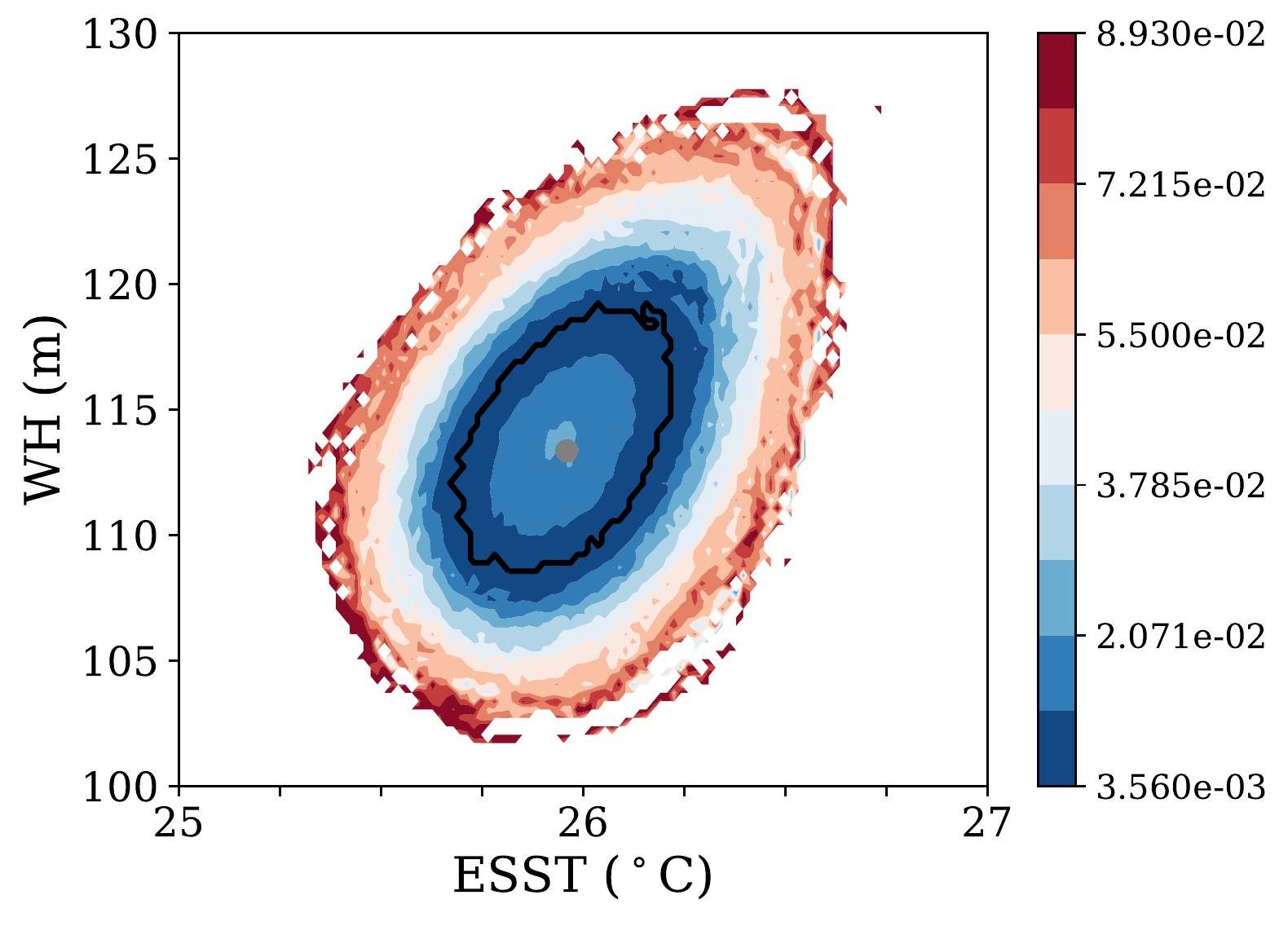}
		\caption{}
	\end{subfigure}\\
	\begin{subfigure}[b]{0.32\textwidth}
		\includegraphics[width=\textwidth]{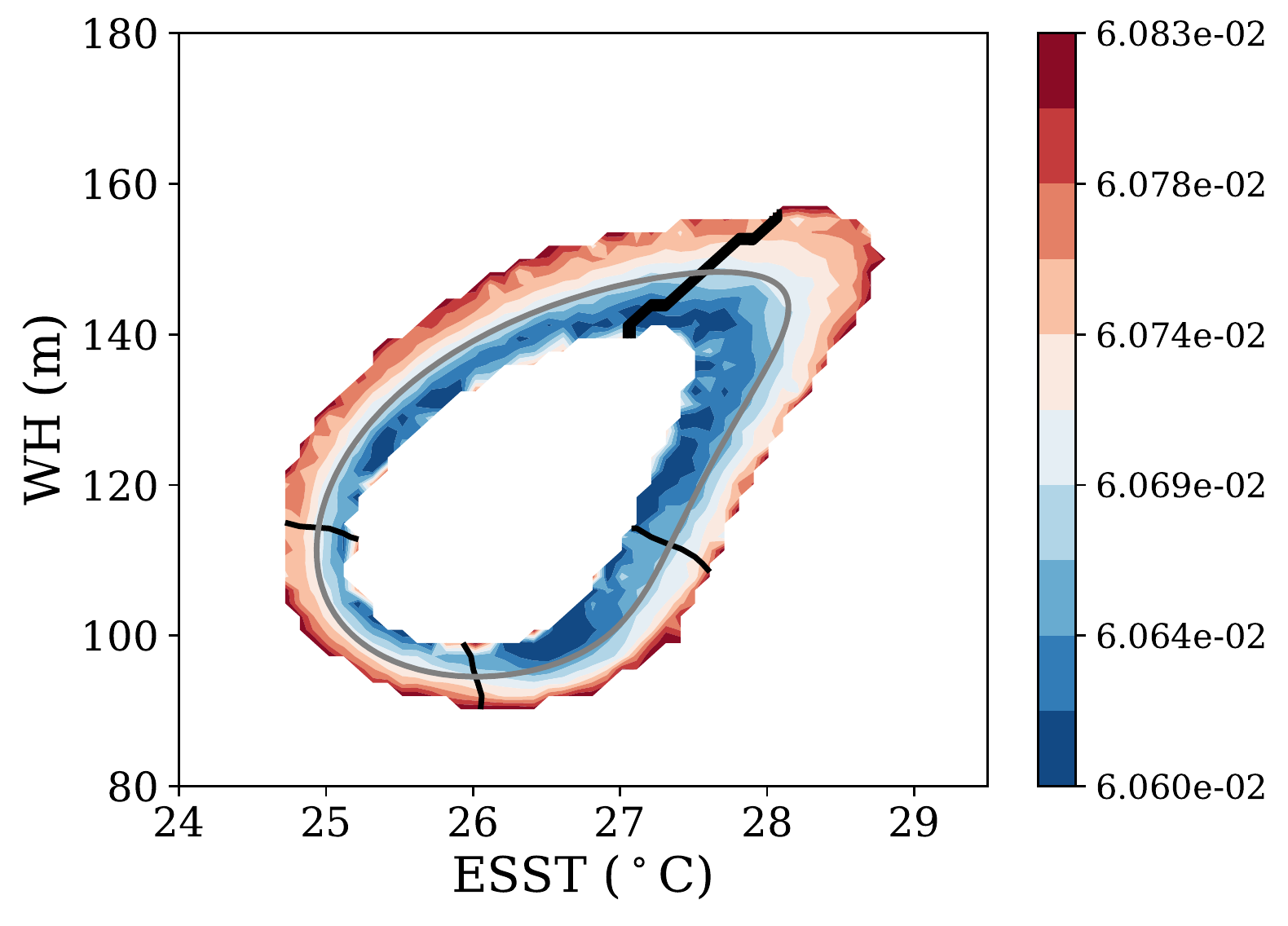}
		\caption{}
	\end{subfigure}
	\begin{subfigure}[b]{0.32\textwidth}
		\includegraphics[width=\textwidth]{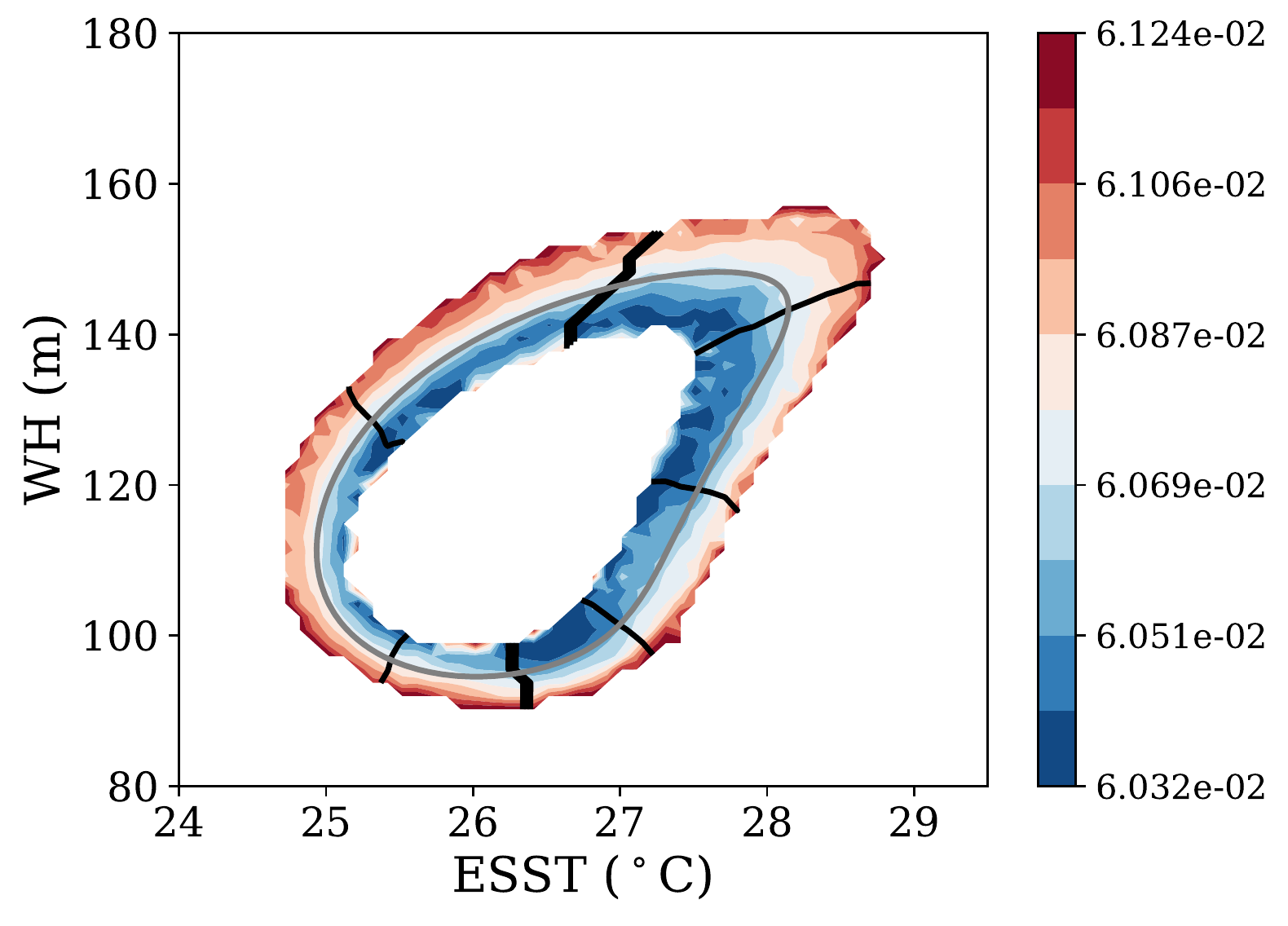}
		\caption{}
	\end{subfigure}
	\begin{subfigure}[b]{0.32\textwidth}
		\includegraphics[width=\textwidth]{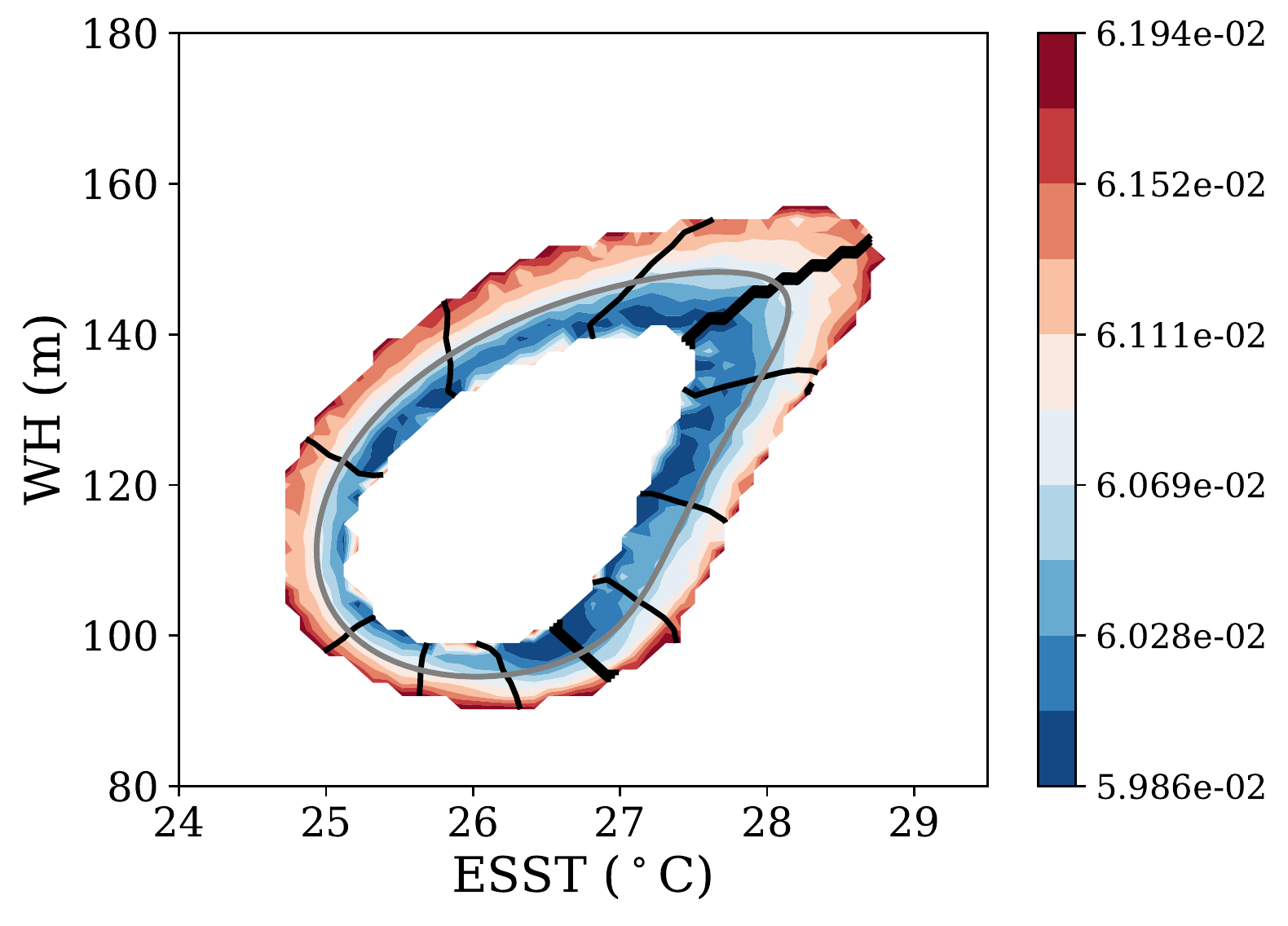}
		\caption{}
	\end{subfigure}\\
	\begin{subfigure}[b]{0.32\textwidth}
		\includegraphics[width=\textwidth]{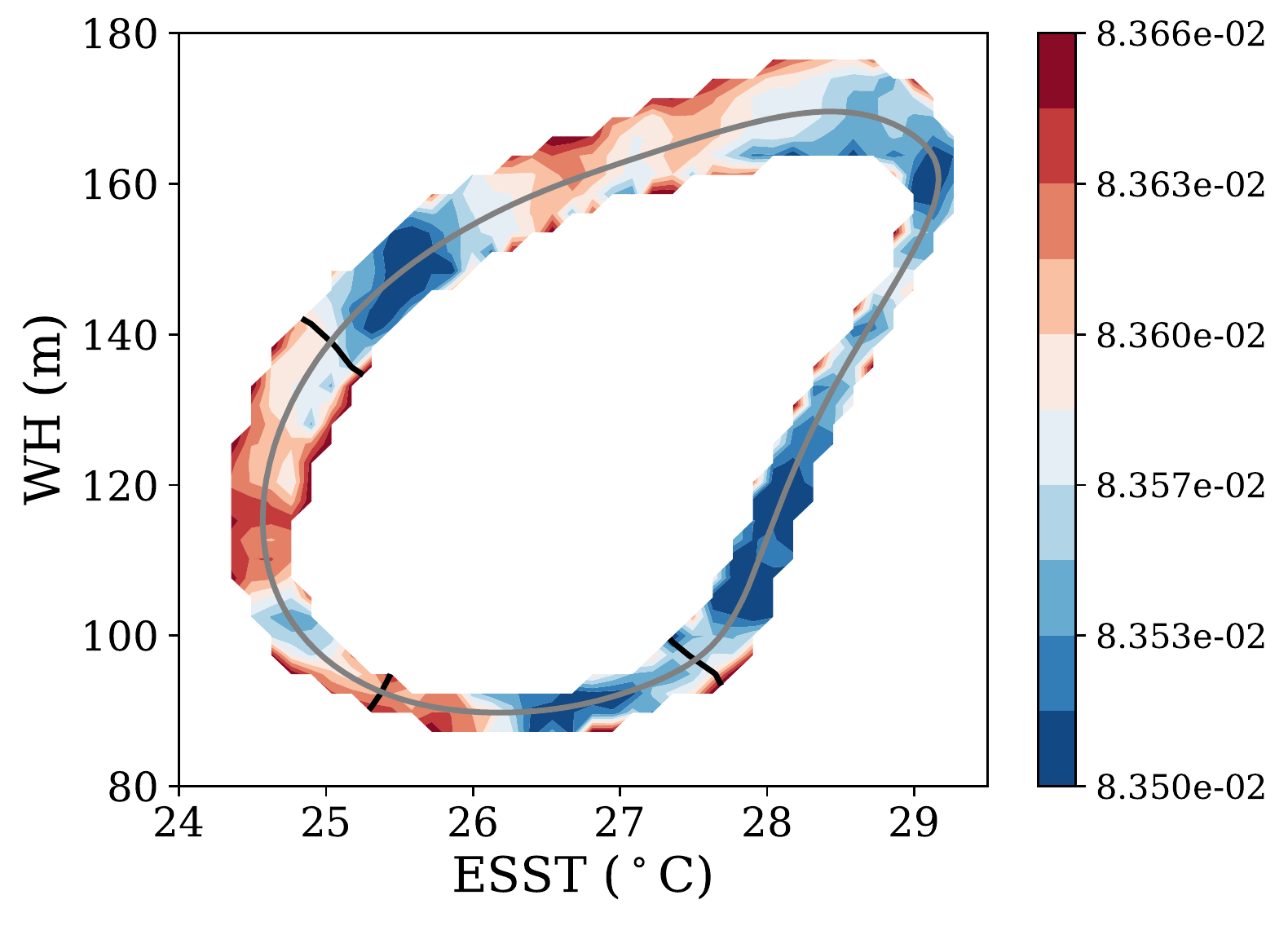}
		\caption{}
	\end{subfigure}
	\begin{subfigure}[b]{0.32\textwidth}
		\includegraphics[width=\textwidth]{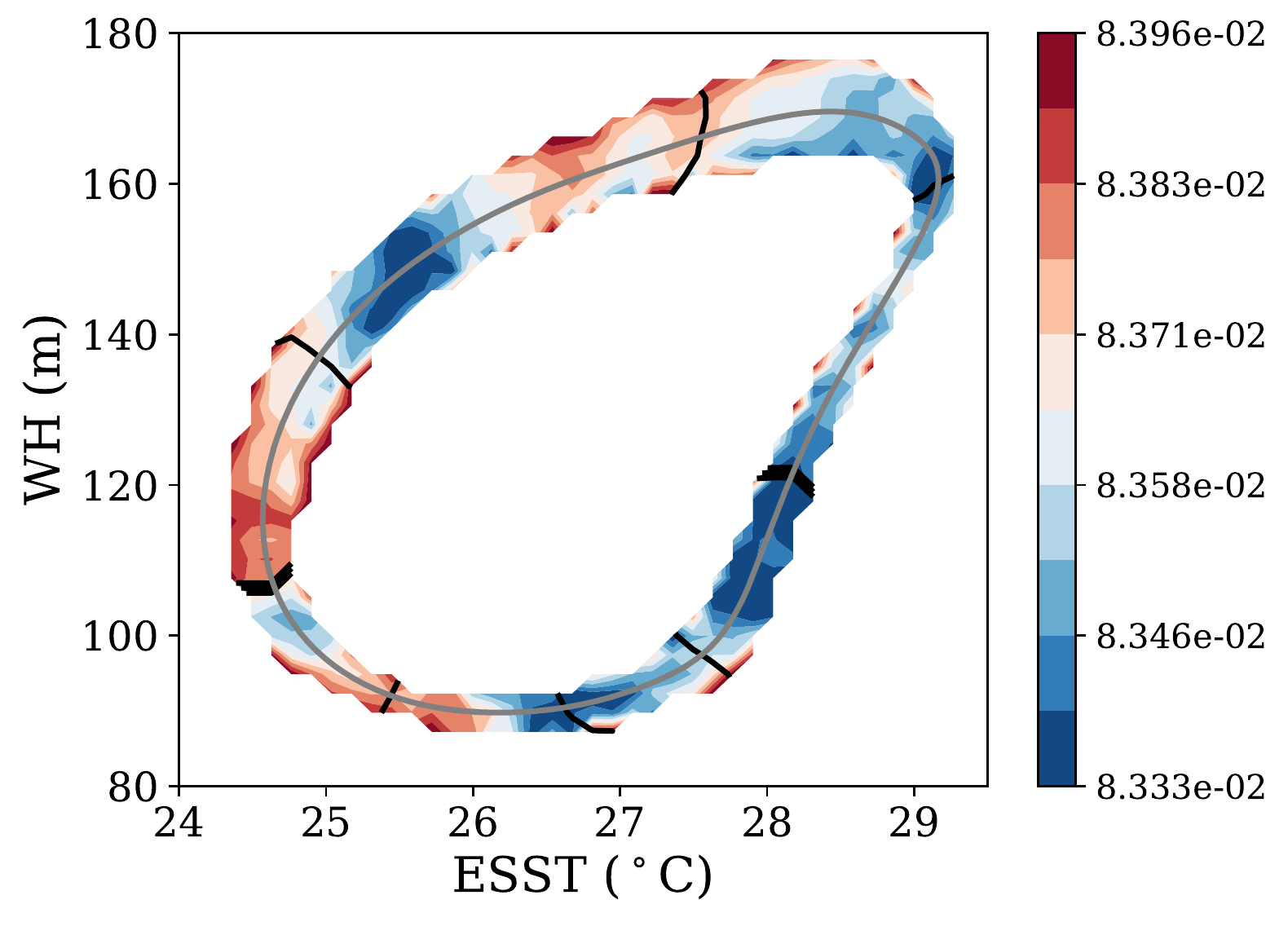}
		\caption{}
	\end{subfigure}
	\begin{subfigure}[b]{0.32\textwidth}
		\includegraphics[width=\textwidth]{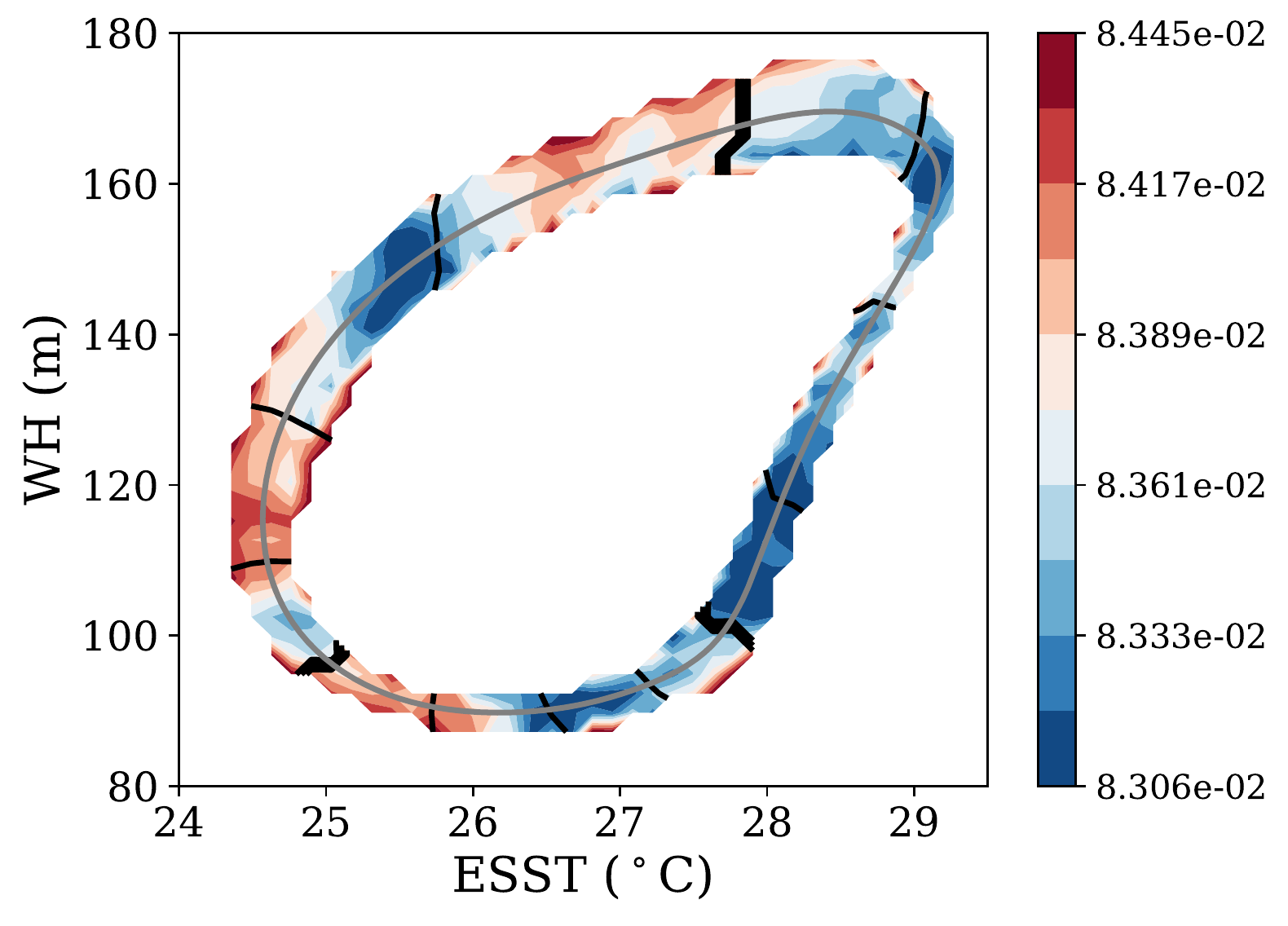}
		\caption{}
	\end{subfigure}\\
	\caption{Argument (line contours with a step of $\pi/2$) and modulus (shaded contours) of the second (left), fourth (center) and sixth (right) reduced backward eigenvectors of the stochastic CZ, for a coupling of
          (a-c) 2.80 (before the deterministic Hopf bifurcation),
          (d-f) 2.85 (right before the bifurcation),
          (g-i) 2.90 (right after the bifurcation),
          (j-l) 2.95 (after the bifurcation).
          The stable stationary point and the limit cycle from the deterministic version of the model and the same values of the coupling are also represented in grey.
          The scale is changed between panels (a--f) and (g--l).
          The closed contour line in panel (b) and (f) indicate a change of sign in the amplitude of the eigenvectors, which are in fact real.}\label{fig:eigvecCZ}
\end{figure}
The eigenvectors associated with the reduced RP resonances are also instructive.
In Fig.~\ref{fig:eigvecCZ}, we represent the 2nd, 4th and 6th eigenvectors (from left to right), for the same values of the coupling as for Fig.~\ref{fig:spectrumCZ} (from top to bottom).
Because of the periodic character of the dynamics, the modulus and the argument of these complex vectors are represented as filled contours and line contours, respectively.
This choice of representation differs from that used in part II of this contribution~\cite{Tantet2017b} and is made to facilitate the lecture of radial variations.
Moreover, the 1st, 3rd and 5th eigenvectors are not represented because the 1st is always constant in space and the 3rd and 5th are the complex conjugates of the 2nd and 4th eigenvectors.

Apart from the 4th eigenvector for $\xi = 2.80$ and the 6th eigenvector for $\xi = 2.85$ (panels b and f) which, being associated with a real eigenvalue, are real, the arguments of the eigenvectors present wavelike patterns with a wavenumber increasing with the rank of the associated eigenvalue (from wavenumber 1 to wavenumber 3).
In agreement with the pseudo-periodic dynamics, these wavenumbers correspond to the rank of the harmonics given by the imaginary part of the associated eigenvalues.
For example, for $\xi = 2.90$, the 2nd eigenvector (panel g) shows a wavenumber 1 pattern and is associated with an eigenvalue with imaginary part $|\Im(\lambda_1)| = \SI{2.00}{rad/y}$, while the 4th (resp.~6th) eigenvector shows a wavenumber 2 (resp.~3) pattern and is associated with an eigenvalue with imaginary part $|\Im(\lambda_3)| = \SI{4.00}{rad/y}$ (resp.~$|\Im(\lambda_3)| = \SI{6.00}{rad/y}$).

Second, one can see that the contour shades are roughly concentric,
so that the modulus' of the eigenvectors increase in the radial direction.
Moreover, the shrinking gap between the contour lines indicates that the gradient in the modulus of the eigenvectors is sharper for eigenvalues of increasing rank.
Also, the modulus tends to flatten as the coupling is increased.

Finally, for large values of $\xi$, the isolines of phase are tilted concentric curves rather that straight lines perpendicular to the limit cycle.

%Indeed, there patterns are notably skewed and concentrated towards the deterministic stationary point
%with a value of $30^\circ C$ for the SSTs, with negative anomalies compared to positive anomalies.
%This can be physically understood by the fact that the temperature of the ocean tends to relax to the
%radiative equilibrium temperature $T_E$ set to $30^\circ$, while, during a La Ni\~na event, the ocean states gets cooler.
%Thus, the reduced mixing eigenvectors represented here bear the signature of the nonlinear dynamics proper to the CZ model
%and which cannot be captured by normal form theory.
%
%This shows that not only the mixing eigenvalues can allow to unravel the nonlinear dynamics in complex models
%but also the patterns of the mixing eigenvectors.
%These important points will be summarized and discussed in greater details in the following section \ref{sec:Conclusion}.

%%%
\section{Interpretation from small-noise expansions}\label{sec:resultsJin}

In this section, we build on small-noise expansions of the RP resonances for stationary points and limit cycles to interpret the spectral results found in the
previous Section~\ref{Sec_mixing_spec_CZ} for the CZ model.
General small-noise expansions have been obtained by~\cite{gaspard2002trace} for hyperbolic systems (i.e.~away from the bifurcation, in the case of a Hopf).
Explicit formulas are derived in the second part of this contribution~\cite{Tantet2017b} for the Hopf normal form subject to an additive white noise, and the associated eigenvectors were also found.
The analytical formulas for the RP resonances associated with a limit cycle forced by noise can be derived, to first orders in the noise level, thanks to the small-noise expansion approach.
To compute the terms of these formulas for a given limit cycle, the Floquet analysis of the limit cycle is needed.
In most cases, this analysis needs to be done numerically, based on the continuation of the limit cycle and on the computation of the corresponding fundamental matrix.
In the case of the CZ model used here, this algorithm has not been implemented. However, the small-noise expansions of the stochastic Hopf bifurcation from~\cite{Tantet2017b} already provide a number of key elements to interpret the results of the previous Section~\ref{Sec_mixing_spec_CZ}.

\subsection{Hopf normal form with noise}

In~\cite{Tantet2017b}, white noise is uniformly added to the system of ordinary differential equations of the Hopf normal form in Cartesian coordinates.
In polar coordinates, the resulting system of SDEs is

\begin{align}
  d r
  &= (\delta r - r^3 + \frac{\epsilon^2}{2 r}) d t  + \epsilon d W_r \\
  d \theta
  &= (\gamma - \beta r^2) d t + \frac{\epsilon}{r} d W_\theta, \label{eq:HopfSDEPolar}
\end{align}
where $W_r$ and $W_\theta$ are two independent Wiener processes
with differentials interpreted in the \^Ito sense.
Here, $\delta$ controls the stability of the asymptotic solutions:
a hyperbolic stationary point at the origin for $\delta < 0$,
a hyperbolic limit cycle of radius $\sqrt{\delta}$ for $\delta > 0$. The parameter $\gamma$ yields a constant contribution to the angular frequency of the periodic solution, while $\beta$ controls the dependence of the angular frequency on the radius.
On the other hand, $\epsilon$ controls the noise intensity.
Small-noise expansions of the RP resonances and their eigenvectors have been derived in~\cite{Tantet2017b} for the SDE~\eqref{eq:HopfSDEPolar} and for $|\delta| > 0$ (i.e.~away from the bifurcation).
We summarize these developments here and use them to interpret the results for the CZ model.

\subsection{Before the bifurcation}

For $\delta < 0$, There exists only one family of RP resonances associated with the hyperbolic stationary point.
They are given to first order in the noise intensity by

\begin{equation}
	\lambda_{ln} = (l + n) \delta + i (n - l) \gamma + \mathcal{O}(\epsilon^2), \quad l, n \in \mathbb{N}.
	\label{eq:eigValFPSub}
\end{equation}
The eigenvalues are thus arranged in a triangular structure to the left of the imaginary axis and generated by linear combinations of the eigenvalues $\delta \pm i \gamma$ of the Jacobian of the deterministic stationary point (i.e~for $\epsilon = 0$).
In fact, to first order, the RP resonances coincide with those of the deterministic case. The associated eigenfunctions are given by

\begin{equation}
	\psi_{ln}(r, \theta) \approx
	\begin{cases}
		e^{i(n-l)\theta} \enskip
		\sqrt{\frac{l!}{n!}}
		\left(\sqrt{-\frac{\delta}{\epsilon^2}} r\right)^{n-l}
		L_l^{n-l}(-\frac{\delta r^2}{\epsilon^2}), \quad	&n \ge l \\
		e^{i(l-n)\theta} \enskip
		\sqrt{\frac{n!}{l!}}
		\left( \sqrt{-\frac{\delta}{\epsilon^2}}r \right)^{l-n}
		L_n^{l-n}(-\frac{\delta r^2}{\epsilon^2}), \quad		&n < l.
	\end{cases}
	\label{eq:complexHermite}
\end{equation}
In other words, the eigenfunctions are the product of a harmonic function in $\theta$ with a Laguerre polynomial $L^\alpha_n$ in $r$.
The order of the harmonic function corresponds to the multiple of the fundamental frequency given by the imaginary part of the eigenvalue (the higher the rank of the harmonic, the larger the wavenumber).
The degree of the polynomial corresponds to the multiple of the stability coefficient $\delta$ giving the real part of the eigenvalue (the further the eigenvalue from the imaginary axis, the higher the degree).

These formulas explain several features of the results (Fig.~\ref{fig:spectrumCZ} and~\ref{fig:eigvecCZ}) obtained in Section~\ref{Sec_mixing_spec_CZ} for the CZ model with $\xi = 2.80, 2.85$ (before the deterministic bifurcation), at least for the leading resonances.
The eigenvalues are arranged in a triangular structure with the imaginary parts (resp.~real parts) of the eigenvalues corresponding to harmonic functions of increasing wave number (resp.~nonlinear functions with an increasing number of changes of sign) in the associated eigenvectors.
The small-noise expansions also explain the approach of the resonances to the imaginary axis with the decreasing stability of the stationary point, as controlled by the parameters $\delta$ (normal form) and $\xi$ (CZ).

\subsection{After the bifurcation}

For $\delta > 0$, two families of eigenvalues coexist.
One is associated with the unstable stationary point at the origin and is again arranged in a triangular structure, but this time with strictly negative real parts.
The other is associated with the deterministic limit cycle and is given by the following small-noise expansions

\begin{equation}
	\lambda_{ln} =
	\begin{cases}
		-\frac{n^2 \epsilon^2 (1 + \beta^2)}{2 \delta} + i n (\gamma - \beta \delta) + o(\epsilon^4),
		\quad &l = 0, \quad n \in \mathbb{Z} \\
		- 2 l \delta + i n (\gamma - \beta \delta) + \mathcal{O}(\epsilon^2), \quad &l \ne 0.
	\end{cases}
	\label{eq:eigValOrbit}
\end{equation}
The eigenvalues are thus arranged in a series of parabolas
such that the real parts of the eigenvalues within a parabola increases with the imaginary part.
The latter are again given by multiples of the fundamental frequency.
The first parabola passes through the eigenvalue 0, which is the only eigenvalue
with zero real part.
Interestingly, the distance of the other eigenvalues from the imaginary axis scales both
with the ratio $\epsilon^2 / (2 \delta)$ and with the coefficient $1 + \beta^2$.
We will see below, that this factor measures the phase diffusion, i.e.~the diffusion in the direction tangent to the limit cycle and responsible for the mixing along it.
This factor relates to the dependence of the angular velocity on the distance to the limit cycle through the coefficient $\beta$.
The other parabolas are separated by a gap given by a multiple of the Floquet exponent $-2 \delta$ of the deterministic limit cycle.

The eigenfunctions associated with the first row of eigenvalues are given by

\begin{equation}
	\psi_{ln} \approx \left(2^k k!\right)^{-\frac{1}{2}} \enskip
	e^{i n \left(\theta - \beta \log{\frac{r}{\sqrt{\delta}}} \right)} \enskip
	H_l\left( \frac{\sqrt{2 \delta}}{\epsilon} \left(r - \sqrt{\delta}\right) \right), \enskip \quad l = 0.
	\label{eq:eigvecHarmonic}
\end{equation}
Once again, the eigenfunctions are the product of a harmonic function
with a polynomial.
The latter are given by Hermite polynomials in $r$, centered at the deterministic limit cycle, and scaled with the variance of the process about the cycle.
The wavenumber is also increasing with the harmonic corresponding to the imaginary part of the eigenvalue and the degree of the polynomial with the distance of the eigenvalue from the imaginary axis.
When $\beta$ is different from zero, the harmonic functions depend on the radius, so that the isolines of phase are tilted.
These lines correspond to the isochrons of the limit cycle (i.e.~the foliation of its stable manifold)~\cite{Mauroy2013}.
In fact, it was shown in~\cite{Tantet2017b} that the noise field has to have a component transverse to the isochrons for phase diffusion to occur.

Several features of the CZ-model results for a coupling of $\xi = 2.90, 2.95$ (after the deterministic bifurcation, Fig.~\ref{fig:spectrumCZ} and~\ref{fig:eigvecCZ}) can be explained from these formulas.
The resonances are indeed found to be arranged in series of parabolic lines of eigenvalues.
The gap between the imaginary parts and the real parts of the resonances correspond to the wave number and the number of changes of sign in the associated eigenvectors.
Moreover, the leading secondary resonances are indeed separated from the imaginary axis by a gap which could be explained by the phase diffusion quantified by~\eqref{eq:eigValOrbit}, and the tilt in the isolines of phase of the eigenvectors is also found in~\eqref{eq:eigvecHarmonic}, for $|\beta > 0|$.
Thus, in the case of a stable limit cycle (after the Hopf bifurcation, here) and for a weak noise, the decay of correlations is primarily governed by the phase diffusion along the orbit. This is different from the case of a stable stationary point (before the Hopf bifurcation, here) for which the rate of decay of correlations is dominated by the stability of this point, and for which the noise level and the phase diffusion only play a secondary role.
On the other hand, the skewness of the eigenvectors of the CZ model is not explained by the small-noise expansion~\eqref{eq:eigvecHarmonic}.
In this case, we expect the Floquet elements in the direction transverse to the limit cycle to play an important role~\cite{Mauroy2013}.

\subsection{General formulas for stationary points and limit cycles}\label{sec:generalFormulas}

To go further in the interpretation of the spectral properties of the CZ model, we now discuss the general formulas of the small-noise expansions derived by~\cite{gaspard2002trace}. Although numerical applications of these formulas require the computation of the Jacobian (resp.~fundamental matrix) of the deterministic stationary point (resp.~limit cycle), which is beyond the scope of this study, they provide further insights on the structure of the resonances of the CZ model.

The RP resonances associated with a hyperbolic stationary point in the presence of white noise are given to first order in the noise intensity by

\begin{align}
	\lambda_{\mathbf{l} \mathbf{n}} = \sum_{\Re(\alpha_i) < 0} l_i \alpha_i
	- \sum_{\Re(\alpha_j) > 0} n_j \alpha_j + \mathcal{O}(\epsilon^2),
	\label{eq:eigValFP}
\end{align}
for $l_i, n_j = 0, 1, 2, \dots$ and where the $\alpha_i, 1 \le i \le n$ are the eigenvalues
of the Jacobian matrix of the stationary point.
In the case of a hyperbolic stable stationary point, $\alpha_i < 0$ for all $i$.
If the first two $\alpha_i$'s form a complex conjugate pair with a small real part,~\eqref{eq:eigValFP} yields a triangular arrangement of RP resonances, as found for the CZ model in Figure~\ref{fig:spectrumCZ} for $\xi = 2.80, 2.85$. The gap between the imaginary (resp.~real) axis and the secondary resonances is thus given by multiples of the real (resp.~imaginary) part of the leading Jacobian eigenvalues.

While the first order of the expansion of the RP resonances associated with a stationary point coincides with the ones of the deterministic problem~\cite{Gaspard1995}, a new phenomena arises in the case of the resonances of a stochastically perturbed hyperbolic limit cycle.
The resonances associated with a hyperbolic limit cycle in the presence of white noise of level $\epsilon$ are given to first orders by

\begin{align}
	\lambda_{\mathbf{l} n} =
	\begin{cases}
	- \Phi n^2 + i n \omega + \mathcal{O}(\epsilon^4)
	\quad &\mathbf{l} = 0, \quad n \in \mathbb{Z} \\
	-\sum_{i = 1}^{n - 1} l_i \nu_i  + i n \omega + \mathcal{O}(\epsilon^2)
	\quad &\mathbf{l} \ne 0, \quad n \in \mathbb{Z},
	\end{cases}
	\label{eq:eigValPO}
\end{align}
where the $\nu_i, 1 \le i \le n-1$ are the Floquet multipliers of the limit cycle with angular frequency $\omega = 2\pi / T$ and period $T$.
The coefficient $\Phi$ of \emph{phase diffusion} is given by

\begin{align}
	\Phi = - \frac{\epsilon^2 \omega^2}{2 T}
	\frac{\langle C(T) \vec{f}, \vec{f} \rangle}{\langle \vec{e}, \vec{f} \rangle}.
	\label{eq:phaseDiffusionCoeff}
\end{align}
It measures the amount of diffusion occurring along the limit cycle.
Its effect is to push eigenvalues --- which in the deterministic case would be on the imaginary axis --- away from the imaginary axis.
The phase diffusion coefficient is thus associated with the asymptotic decay of correlations, or mixing in state space.
In~\eqref{eq:phaseDiffusionCoeff}, $\vec{f}$ is left eigenvector
of the Floquet matrix in the direction of the flow (i.e.~tangent to the limit cycle),
$\vec{e}$ is the corresponding right eigenvector and the matrix $C(T)$ is given by the matrix

\begin{align}
	C(t) = \int_0^t M(t) M(-s) D(s) M(-s)^* M(t)^* ds,
	\label{eq:correlationMatrix}
\end{align}
where the integral is taken over a full cycle on the limit cycle, $M(t)$ is the map at time $t$ of the fundamental matrix of the limit cycle and $D(s)$ is the evaluation at time $t$ of the diffusion matrix of the stochastic system\footnote{The matrix $C(t)$ is in fact the \emph{covariance matrix} of a periodic Ornstein-Uhlenbeck process with a drift given by the Jacobian matrix $A(t)$ generating the fundamental matrix $M(t)$ and with the diffusion matrix $D(t)$, both evaluated along the limit cycle.}.

In the case of the stochastic Hopf bifurcation above, one has for $\delta > 0$ that $\omega = \gamma - \beta \delta$ and $\nu_1 = 2 \delta$.
The phase diffusion coefficient is given by $\Phi = \epsilon^2 (1 + \beta^2) / (2 \delta)$.
It is larger the stronger the noise is with respect to the contraction to the limit cycle, but it is also affected by the twist factor $\beta$.
By twisting the isochrons, as apparent in Figure~\ref{fig:eigvecCZ}, the amount of noise in the radial direction which is transferred to the azimuthal direction by the deterministic dynamics is increased~\cite{Tantet2017b}.
In the case of the CZ model, for $\xi = 2.90, 2.95$, we thus expect the gap observed between the leading secondary eigenvalues and the imaginary axis (Fig.~\ref{fig:spectrumCZ}) to be in fact controlled by the phase diffusion coefficient~\eqref{eq:phaseDiffusionCoeff} times the rank squared of the harmonic associated with the resonance.

Finally, as the bifurcation is neared, e.g.~for $\xi = 2.85, 2.9$, higher order terms in the noise level than the ones given here may be necessary to explain the spectral properties of the Cane-Zebiak model.

\section{Conclusion}\label{sec:conclusion}

The response of a low-frequency mode of climate variability, El Ni\~no-Southern Oscillation, to stochastic forcing is studied in a model of intermediate complexity, the fully-coupled Cane-Zebiak model, from the spectral analysis of the Markov semigroup.
In this mathematical framework, resonances in the power spectral density and the associated slow time-decay of correlations are characterized by eigenvalues and eigenfunctions of Markov operators.
For such high-dimensional stochastic systems, extraction of partial information about these operators is made possible by the reduction approach of~\cite{Chekroun2017a}, whereby transition matrices corresponding to the projection of the Markov operators on a low-dimensional space are estimated.
Our interpretations of the numerical results are supported by the Hopf normal form small-noise expansions derived in~\cite{Tantet2017b}.

Without noise, the Cane-Zebiak model undergoes a Hopf bifurcation as the coupling between the ocean and the atmosphere is increased. The addition of noise in the wind stress leads to noise-induced oscillations, even before the critical value of the deterministic bifurcation~\cite{Roulston2000}.
Our results show that, as opposed to the deterministic case for which a sharp bifurcation occurs at a specific value of the coupling, the noise is responsible for a smooth transition of the RP resonances as the coupling is increased.
Thus, it is not possible to define a precise value at which the stochastic Hopf bifurcation occurs.

For low values of the coupling strength and in the presence of noise, our estimations reveal the presence of a triangular structure of eigenvalues, or reduced RP resonances, in the complex plane.
The leading complex conjugate pair of eigenvalues is associated with power spectral peaks reminiscent of that found from ENSO indices.
The leading eigenvectors are centered about the deterministic stationary point, with a phase resembling a harmonic function with a wave number corresponding to the imaginary part of the eigenvalue, and a modulus with an increasing number of sign changes for eigenvalues further from the imaginary axis.
These results are well predicted by small-noise expansions of the RP spectrum for hyperbolic stationary points with complex Jacobian eigenvalues.

As the critical value at which the deterministic bifurcation occurs is approached, the spectral gap between the leading eigenvalues and the imaginary axis shrinks. This explains (i) the sharpening of the peaks in the power spectral density associated with the fundamental frequency of the system and its harmonics, (ii) the corresponding slowing down of the decay of correlations and (iii) the triggering of oscillations by the noise.
As the bifurcation is passed, the leading eigenvalues are rearranged on successive parabolic curves.
In the absence of noise, the leading eigenvalues would lie on the imaginary axis, with their imaginary parts corresponding to the harmonics of the limit cycle.
Instead, in the presence of noise, mixing in state space from phase diffusion is characterized by a the spectral gap between the leading eigenvalues and the imaginary axis (with the exception of the zero eigenvalue). To first order in the noise, the coefficient of phase diffusion measuring this gap can in fact be computed from the Floquet analysis of the deterministic limit cycle~\cite{gaspard2002trace,Tantet2017b}.
The leading eigenvectors now spread about the deterministic cycle.
Their amplitude is constant and there phase resembles harmonic functions of increasing wave number.
Contrary to the case of the stationary point, the isolines of phase, or isochrons, are twisted.
In particular, it is shown in~\cite{Tantet2017b} that the noise field has to have a component transverse to the isochrons for phase diffusion to occur.
For instance, the stronger the twist the more efficient the noise transfer from directions orthogonal to the limit cycle to the direction tangent to it.

From a physical point of view, the study of the RP resonances allows one to better understand how fast atmospheric perturbations may be responsible for the irregularity in the recurrence of ENSO events.
To first order, the sensitivity of ENSO --- measured as the amount of phase diffusion --- would not only be determined by the intensity of the perturbation, but also by the degree of the interaction of the stochastic forcing with the nonlinear ENSO dynamics.
In particular, if the ENSO evolution is found to depend on the amplitude of ENSO events, then its future evolution may be particularly sensitive to atmospheric perturbations and potentially less predictable.

In this study, the ENSO model considered is autonomous, so that the impact of the seasonal cycle is not considered.
Yet, it is known that the seasonal cycle is responsible for the phase locking of ENSO events to peak around December/January and for making predictions before Spring particularly difficult~\cite{Neelin1998a}.
To better understand the joint effect of the seasonal cycle and stochastic forcing on ENSO, the ergodic theoretic framework presented in this study needs to be extended.

%%%%%%%%%%%%%%%%%%%%%%%%%%%%%%%%%%%%%%%%%%%%%%%%%%%%%%%%%%%%%%%%%%%%%
% ACKNOWLEDGMENTS
%%%%%%%%%%%%%%%%%%%%%%%%%%%%%%%%%%%%%%%%%%%%%%%%%%%%%%%%%%%%%%%%%%%%%
%
%\acknowledgments
% Start acknowledgments here.
\begin{acknowledgements}
% The authors would like to thank the reviewers for there careful reading and comments.
This work has been partially supported by the Office of Naval Research (ONR) Multidisciplinary University Research Initiative (MURI) grant N00014-12-1-0911 and N00014-16-1-2073 (MDC), by the National Science Foundation grant DMS-1616981(MDC) and AGS-1540518 (JDN),  by  the LINC project (No.~289447) funded by EC's Marie-Curie ITN (FP7-PEOPLE-2011-ITN) program (AT and HD) and by the Utrecht University Center for Water, Climate and Ecosystems (AT). 
\end{acknowledgements}

\bibliographystyle{amsalpha}
\bibliography{atantet}

\newcommand{\etalchar}[1]{$^{#1}$}
\providecommand{\bysame}{\leavevmode\hbox to3em{\hrulefill}\thinspace}
\providecommand{\MR}{\relax\ifhmode\unskip\space\fi MR }
\providecommand{\MRhref}[2]{%
  \href{http://www.ams.org/mathscinet-getitem?mr=#1}{#2}
}
\providecommand{\href}[2]{#2}
\begin{thebibliography}{FPGET07}

\bibitem[ALS13]{avram_spectral_2013}
F.~Avram, N.~N. Leonenko, and N.~Suvak, \emph{On {{Spectral Analysis}} of
  {{Heavy}}-{{Tailed Kolmogorov}}-{{Pearson Diffusions}}}, no.~19, 249--298.

\bibitem[Bag14]{Bagheri2014}
S.~Bagheri, \emph{{Effects of weak noise on oscillating flows: Linking quality
  factor, Floquet modes, and Koopman spectrum}}, Phys. Fluids \textbf{26}
  (2014), no.~9, 094104.

\bibitem[BKJ15]{Bittracher2015}
A.~Bittracher, P.~Koltai, and O.~Junge, \emph{Pseudogenerators of spatial
  transfer operators}, SIAM J. Appl. Dyn. Syst. \textbf{14} (2015), no.~3,
  1478--1517.

\bibitem[BL07]{Butterley_Liverani}
O.~Butterley and C.~Liverani, \emph{{Smooth Anosov flows: Correlation spectra
  and stability}}, {Journal of Modern Dynamics} \textbf{1} (2007), no.~2,
  301--322.

\bibitem[CCH{\etalchar{+}}16]{chen2016diversity}
C.~Chen, M.~A. Cane, N.~Henderson, D.~E. Lee, D.~Chapman, D.~Kondrashov, and
  M.~D. Chekroun, \emph{{Diversity, nonlinearity, seasonality, and memory
  effect in ENSO simulation and prediction using empirical model reduction}},
  Journal of Climate \textbf{29} (2016), no.~5, 1809--1830.

\bibitem[CCHT19]{cao2019mathematical}
Y.~Cao, M.~D. Chekroun, A.~Huang, and R.~Temam, \emph{Mathematical analysis of
  the jin-neelin model of el ni{\~{n}}o-southern oscillation}, Chinese Annals
  of Mathematics, Series B \textbf{40} (2019), no.~1, 1--38.

\bibitem[Cer01]{cerrai2001second}
S.~Cerrai, \emph{{Second-Order PDE's in Finite and Infinite Dimension: A
  Probabilistic Approach}}, vol. 1762, Springer Science \& Business Media,
  2001.

\bibitem[CNK{\etalchar{+}}14]{Chekroun2014}
M.~D. Chekroun, J.D. Neelin, D.~Kondrashov, J.C. McWilliams, and M.~Ghil,
  \emph{{Rough parameter dependence in climate models: The role of
  Ruelle-Pollicott resonances}}, {Proc. Natl. Acad. Sci} \textbf{111} (2014),
  no.~5, 1684--1690.

\bibitem[CTND19]{Chekroun2017a}
M.~D. Chekroun, A.~Tantet, J.~D. Neelin, and H.~A. Dijkstra,
  \emph{{Ruelle-Pollicott Resonances of Stochastic Systems in Reduced State
  Space. Part I: Theory}}, Submitted to J. Stat. Phys. (2019).

\bibitem[CVE06a]{crommelin2006reconstruction}
D.~Crommelin and E.~Vanden-Eijnden, \emph{Reconstruction of diffusions using
  spectral data from time series}, Communications in Mathematical Sciences
  \textbf{4} (2006), no.~3, 651--668.

\bibitem[CVE06b]{crommelin2006fitting}
DT~Crommelin and Eric Vanden-Eijnden, \emph{Fitting time series by
  continuous-time markov chains: A quadratic programming approach}, Journal of
  Computational Physics \textbf{217} (2006), no.~2, 782--805.

\bibitem[DAXP10]{Deser2010}
Clara Deser, Michael~a. Alexander, Shang-Ping Xie, and Adam~S. Phillips,
  \emph{{Sea Surface Temperature Variability: Patterns and Mechanisms}}, Ann.
  Rev. Mar. Sci. \textbf{2} (2010), no.~1, 115--143.

\bibitem[DDJS99]{Deuflhard1999}
Peter Deuflhard, Michael Dellnitz, Oliver Junge, and Christof Sch{\"{u}}tte,
  \emph{{Computation of Essential Molecular Dynamics by Subdivision
  Techniques}}, Comput. Mol. Dyn. Challenges, Methods, Ideas (Peter Deuflhard,
  Jan Hermans, Benedict Leimkuhler, Alan~E. Mark, Sebastian Reich, and
  Robert~D. Skeel, eds.), vol.~45, Springer, Berlin, 1999, pp.~98--115.

\bibitem[DFH{\etalchar{+}}09]{Dellnitz2009}
Michael Dellnitz, Gary Froyland, C.~Horenkamp, Kathrin Padberg-Gehle, and
  Alexander {Sen Gupta}, \emph{{Seasonal variability of the subpolar gyres in
  the Southern Ocean: a numerical investigation based on transfer operators}},
  Nonlinear Process. Geophys. \textbf{16} (2009), no.~6, 655--663.

\bibitem[Dij13]{Dijkstra2013}
Henk~A. Dijkstra, \emph{{Nonlinear Climate Dynamics}}, Cambridge University
  Press, Cambridge, 2013.

\bibitem[DJ97]{Dellnitz1997a}
Michael Dellnitz and Oliver Junge, \emph{{Almost invariant sets in Chua's
  circuit}}, Int. J. Bifurc. Chaos \textbf{7} (1997), no.~11, 2475--2485.

\bibitem[DJ99]{Dellnitz1999}
\bysame, \emph{{On the Approximation of Complicated Dynamical Behavior}}, SIAM
  J. Numer. Anal. \textbf{36} (1999), no.~2, 491--515.

\bibitem[DZ15]{dyatlov2015stochastic}
S.~Dyatlov and M.~Zworski, \emph{{Stochastic stability of Pollicott--Ruelle
  resonances}}, Nonlinearity \textbf{28} (2015), no.~10, 3511.

\bibitem[EN01]{Engel2001}
Klaus-Jochen Engel and Rainer Nagel, \emph{{One-parameter semigroups for linear
  evolution equations}}, Springer, New York, 2001.

\bibitem[FGH13]{Froyland2013}
Gary Froyland, Georg~A. Gottwald, and Andy Hammerlindl, \emph{{A computational
  method to extract macroscopic variables and their dynamics in multiscale
  systems}}, arXiv Prepr. arXiv1310.8001 \textbf{13} (2013), no.~4, 1--30.

\bibitem[FGH14]{froyland2014computational}
G.~Froyland, G.~A. Gottwald, and A.~Hammerlindl, \emph{A computational method
  to extract macroscopic variables and their dynamics in multiscale systems},
  SIAM J. Appl. Dyn. Syst. \textbf{13} (2014), no.~4, 1816--1846.

\bibitem[FGP10]{flandoli2010flow}
F.~Flandoli, M.~Gubinelli, and E.~Priola, \emph{{Flow of diffeomorphisms for
  SDEs with unbounded Holder continuous drift}}, Bulletin des sciences
  mathematiques \textbf{134} (2010), no.~4, 405--422.

\bibitem[FPG09]{Froyland2009}
Gary Froyland and Kathrin Padberg-Gehle, \emph{{Almost-invariant sets and
  invariant manifolds - Connecting probabilistic and geometric descriptions of
  coherent structures in flows}}, Phys. D Nonlinear Phenom. \textbf{238}
  (2009), no.~16, 1507--1523.

\bibitem[FPGET07]{Froyland2007}
Gary Froyland, Kathrin Padberg-Gehle, Matthew England, and Anne Treguier,
  \emph{{Detection of Coherent Oceanic Structures via Transfer Operators}},
  Phys. Rev. Lett. \textbf{98} (2007), no.~22, 224503.

\bibitem[Fro97]{froyland1997computer}
G.~Froyland, \emph{Computer-assisted bounds for the rate of decay of
  correlations}, Communications in Mathematical Physics \textbf{189} (1997),
  no.~1, 237--257.

\bibitem[FSvS14]{Froyland2014b}
Gary Froyland, Robyn~M. Stuart, and Erik van Sebille, \emph{{How well-connected
  is the surface of the global ocean?}}, Chaos An Interdiscip. J. Nonlinear
  Sci. \textbf{24} (2014), no.~3, 033126.

\bibitem[Gas02]{gaspard2002trace}
P.~Gaspard, \emph{Trace formula for noisy flows}, {Journal of Statistical
  Physics} \textbf{106} (2002), no.~1-2, 57--96.

\bibitem[GNPT95]{Gaspard1995}
Pierre Gaspard, G.~Nicolis, A.~Provata, and S.~Tasaki, \emph{{Spectral
  signature of the Pitchfork bifurcation: Liouville equation approach}}, Phys.
  Rev. E \textbf{51} (1995), no.~1, 74--94.

\bibitem[GO81]{Goldenberg1981}
Stanley~B. Goldenberg and James~J. O'Brien, \emph{{Time and space variability
  of tropical Pacific wind stress}}, Mon. Weather Rev. \textbf{109} (1981),
  1190--1207.

\bibitem[Jin96]{Jin1996b}
Fei-Fei Jin, \emph{{Tropical Ocean-Atmosphere Interaction, the Pacific Cold
  Tongue, and the El Ni{\~{n}}o-Southern Oscillation}}, Science (80-. ).
  \textbf{274} (1996), 76.

\bibitem[JN93a]{Jin_al93_part1}
F.-F. Jin and J.~D. Neelin, \emph{{Modes of interannual tropical
  ocean-atmosphere interaction-A unified view. Part I: Numerical results}},
  Journal of the atmospheric sciences \textbf{50} (1993), no.~21, 3477--3503.

\bibitem[JN93b]{Jin_al93_part3}
\bysame, \emph{{Modes of interannual tropical ocean-atmosphere interaction-A
  unified view. Part III: Analytical results in fully coupled cases}}, Journal
  of the atmospheric sciences \textbf{50} (1993), no.~21, 3523--3540.

\bibitem[KKS15]{Klus2015a}
Stefan Klus, P{\'{e}}ter Koltai, and Christof Sch{\"{u}}tte, \emph{{On the
  numerical approximation of the Perron–Frobenius and Koopman operator}},
  arXiv (2015), 1--19.

\bibitem[Kol10]{Koltai2010}
P{\'{e}}ter Koltai, \emph{{Efficient approximation methods for the global
  long-term behavior of dynamical systems - Theory, algorithms and examples}},
  Ph.D. thesis, Technische Universit{\"{a}}t at M{\"{u}}nchen, 2010, p.~162.

\bibitem[LSY97]{Lehoucq1997}
Richard~B. Lehoucq, D.~C. Sorensen, and C.~Yang, \emph{{ARPACK Users' Guide:
  Solution of Large Scale Eigenvalue Problems with Implicitly Restarted Arnoldi
  Methods}}, 1997, pp.~xv + 137.

\bibitem[MM16]{Mauroy2016}
A.~Mauroy and I.~Mezi{\'{c}}, \emph{{Global Stability Analysis Using the
  Eigenfunctions of the Koopman Operator}}, IEEE Transactions on Automatic
  Control \textbf{61} (2016), no.~11, 3356--3369.

\bibitem[MMM13]{Mauroy2013}
A.~Mauroy, I.~Mezi{\'{c}}, and J.~Moehlis, \emph{{Isostables, isochrons, and
  Koopman spectrum for the action-angle representation of stable fixed point
  dynamics}}, Physica D \textbf{261} (2013), 19--30.

\bibitem[ND95]{Neelin1995}
J.~David Neelin and Henk~A. Dijkstra, \emph{{Ocean-Atmosphere Interaction and
  the Tropical Climatology Part I: The Dangers of Flux Correction}}, J. Clim.
  \textbf{8} (1995), no.~5, 1325--1342.

\bibitem[NHJ{\etalchar{+}}98]{Neelin1998a}
D.~S. Neelin, J. D.~and~Battisti, A.~C. Hirst, F.-F. Jin, Y.~Wakata,
  T.~Yamagata, and S.~E. Zebiak, \emph{{ENSO theory}}, J. Geophys. Res. Ocean.
  \textbf{103} (1998), no.~C7, 14261--14290.

\bibitem[NJ93]{Jin_al93_part2}
J.~D. Neelin and F.-F. Jin, \emph{{Modes of interannual tropical
  ocean-atmosphere interaction-a unified view. Part II: Analytical results in
  the weak-coupling limit}}, Journal of the atmospheric sciences \textbf{50}
  (1993), no.~21, 3504--3522.

\bibitem[Pav14]{Pavliotis2014}
Grigorios~A. Pavliotis, \emph{{Stochastic Processes and Applications}},
  Springer, New York, 2014.

\bibitem[Pol86]{pollicott1986meromorphic}
M.~Pollicott, \emph{Meromorphic extensions of generalised zeta functions},
  {Inventiones Mathematicae} \textbf{85} (1986), no.~1, 147--164.

\bibitem[RMB{\etalchar{+}}09]{rowley2009spectral}
C.~W. Rowley, I.~Mezi{\'c}, S.~Bagheri, P.~Schlatter, and D.~S. Henningson,
  \emph{Spectral analysis of nonlinear flows}, Journal of Fluid Mechanics
  \textbf{641} (2009), 115--127.

\bibitem[RN00]{Roulston2000}
Mark~S. Roulston and J.~David Neelin, \emph{{The response of an ENSO Model to
  climate noise, weather noise and intraseasonal forcing}}, Geophys. Res. Lett.
  \textbf{27} (2000), no.~22, 3723--3726.

\bibitem[Rue86]{ruelle1986locating}
D.~Ruelle, \emph{Locating resonances for axiom a dynamical systems}, {Journal
  of Statistical Physics} \textbf{44} (1986), no.~3-4, 281--292.

\bibitem[Sch10]{schmid2010dynamic}
P.~J Schmid, \emph{Dynamic mode decomposition of numerical and experimental
  data}, Journal of Fluid Mechanics \textbf{656} (2010), 5--28.

\bibitem[SFHD99a]{schutte1999direct}
Ch. Sch{\"u}tte, A.~Fischer, W.~Huisinga, and P.~Deuflhard, \emph{A direct
  approach to conformational dynamics based on hybrid monte carlo}, Journal of
  Computational Physics \textbf{151} (1999), no.~1, 146--168.

\bibitem[SFHD99b]{Schutte1999}
Christof Sch{\"{u}}tte, A~Fischer, Wilhelm Huisinga, and P~Deuflhard, \emph{{A
  Direct Approach to Conformational Dynamics Based on Hybrid Monte Carlo}}, J.
  Comput. Phys. \textbf{151} (1999), no.~1, 146--168.

\bibitem[SR96]{Smith1996}
T.~M. Smith and RW~Reynolds, \emph{{Reconstruction of Historical Sea Surface
  Temperatures Using Empirical Orthogonal Functions}}, J. Clim. \textbf{9}
  (1996), 1403--1420.

\bibitem[TCND19]{Tantet2017b}
A.~Tantet, M.~D. Chekroun, J.~D. Neelin, and H.~A. Dijkstra,
  \emph{{Ruelle-Pollicott Resonances of Stochastic Systems in Reduced State
  Space. Part II: Stochastic Hopf Bifurcation}}, Submitted to J. Stat. Phys.
  (2019).

\bibitem[TLD18]{tantet_resonances_2018}
A.~Tantet, V.~Lucarini, and H.~A. Dijkstra, \emph{Resonances in a {{Chaotic
  Attractor Crisis}} of the {{Lorenz Flow}}}, J. Stat. Phys. \textbf{170}
  (2018), no.~3, 584--616.

\bibitem[TLLD18]{tantet_crisis_2018}
A.~Tantet, V.~Lucarini, F.~Lunkeit, and H.A. Dijkstra, \emph{Crisis of the
  {{Chaotic Attractor}} of a {{Climate Model}}: {{A Transfer Operator
  Approach}}}, Nonlinearity \textbf{31} (2018), no.~5, 2221.

\bibitem[TRL{\etalchar{+}}13]{Tu2013}
Jonathan~H. Tu, Clarence~W. Rowley, Dirk~M. Luchtenburg, Steven~L. Brunton, and
  J.~Nathan Kutz, \emph{{On Dynamic Mode Decomposition: Theory and
  Applications}}, 1--30.

\bibitem[TvdBD15]{tantet_early_2015}
A.~Tantet, F.~R. van~der Burgt, and H.~A. Dijkstra, \emph{An early warning
  indicator for atmospheric blocking events using transfer operators}, Chaos
  \textbf{25} (2015), no.~3, 036406.

\bibitem[Ula64]{ulam1964collection}
S.~M. Ulam, \emph{{Problems in Modern Mathematics}}, science ed., Wiley, New
  York, 1964.

\bibitem[vdVDJ00]{Vaart2000}
P~van~der Vaart, Henk~A. Dijkstra, and Fei-Fei Jin, \emph{{The Pacific Cold
  Tongue and the ENSO Mode : A Unified Theory within the Zebiak - Cane Model}},
  J. Atmos. Sci. \textbf{57} (2000), 967--988.

\bibitem[VM08]{Vaidya2008}
Umesh Vaidya and Prashant~G. Mehta, \emph{{Lyapunov measure for almost
  everywhere stability}}, IEEE Trans. Automat. Contr. \textbf{53} (2008),
  no.~1, 307--323.

\bibitem[vSEFS12]{Sebille2012}
Erik van Sebille, Matthew~H England, Gary Froyland, and Erik~Van Sebille,
  \emph{{Origin, dynamics and evolution of ocean garbage patches from observed
  surface drifters}}, Environ. Res. Lett. \textbf{7} (2012), no.~4, 044040.

\bibitem[vSZ99]{VonStorch1999b}
Hans von Storch and F.~Zwiers, \emph{{Stastistical Analysis in Climate
  Research}}, Cambridge University Press, Cambridge, 1999.

\bibitem[Wie85]{wiesenfeld_noisy_1985}
Kurt Wiesenfeld, \emph{Noisy precursors of nonlinear instabilities}, Journal of
  Statistical Physics \textbf{38} (1985), no.~5, 1071--1097 (en).

\bibitem[WK82]{wiesenfeld_effect_1982}
K.~A. Wiesenfeld and E.~Knobloch, \emph{Effect of noise on the dynamics of a
  nonlinear oscillator}, Physical Review A \textbf{26} (1982), no.~5,
  2946--2953.

\bibitem[WKR15]{Williams2015}
Matthew~O. Williams, Ioannis~G. Kevrekidis, and Clarence~W. Rowley, \emph{{A
  Data-Driven Approximation of the Koopman Operator: Extending Dynamic Mode
  Decomposition}}, J. Nonlinear Sci. \textbf{25} (2015), no.~6, 1307--1346.

\bibitem[ZC87]{Zebiak1987a}
Stephen~E. Zebiak and Mark~A. Cane, \emph{{A model of El Nino-Southern
  Oscillation}}, Mon. Weather Rev. \textbf{115} (1987), no.~31, 2262--2278.

\bibitem[Zwo17]{Zworski2017}
Maciej Zworski, \emph{Mathematical study of scattering resonances}, Bulletin of
  Mathematical Sciences \textbf{7} (2017), no.~1, 1--85.

\end{thebibliography}

\end{document}